\documentclass{amsart}
\usepackage{amssymb,amsmath}
\usepackage{array}
\usepackage{url}

\theoremstyle{plain}
\newtheorem{theorem}{Theorem}
\newtheorem{prop}[theorem]{Proposition}
\newtheorem{lem}[theorem]{Lemma}
\newtheorem{cor}[theorem]{Corollary}
\newtheorem{conj}[theorem]{Conjecture}
\newtheorem{conjcor}[theorem]{Corollary of Conjecture}

\theoremstyle{definition}
\newtheorem{defn}{Definition}
\newtheorem{exa}{Example}

\theoremstyle{remark}
\newtheorem{rem}{Remark}

\newcommand{\AAA}{\mathbb{A}}
\newcommand{\ZZZ}{\mathbb{Z}}
\newcommand{\NNN}{\mathbb{N}}
\newcommand{\PPP}{\mathbb{P}}
\newcommand{\FFF}{\mathbb{F}}
\newcommand{\QQQ}{\mathbb{Q}}
\newcommand{\CCC}{\mathbb{C}}

\newcommand{\LLL}{\mathcal L}

\newcommand{\MMMM}{\mathcal M}
\newcommand{\FFFF}{\mathcal F}
\newcommand{\CCCC}{\mathcal C}

\newcommand{\eps}{\varepsilon}
\renewcommand{\phi}{\varphi}

\newcommand{\opname}{\operatorname}
\newcommand{\sep}{\opname{sep}}
\newcommand{\Hom}{\opname{Hom}}
\newcommand{\Morph}{\opname{Mor}}
\newcommand{\Ind}{\opname{Ind}}
\newcommand{\soc}{\opname{soc}}
\newcommand{\socAb}{\opname{soc^{ab}}}
\newcommand{\rank}{\opname{rk}}
\newcommand{\Aut}{\opname{Aut}}
\newcommand{\PGL}{\opname{PGL}}

\newcommand{\cdim}{\opname{covdim}}
\newcommand{\edim}{\opname{edim}}
\newcommand{\rdim}{\opname{rdim}}
\newcommand{\candim}{\opname{cd}}
\newcommand{\Chr}{\opname{char}}
\newcommand{\eval}{\opname{ev}}
\newcommand{\image}{\opname{im}}

\newcommand{\identity}{\opname{Id}}
\newcommand{\ord}{\opname{ord}}
\newcommand{\tdeg}{\opname{tdeg}}
\newcommand{\grk}{\rank_{\ZZZ \!G}}

\newcommand{\GL}{\opname{GL}}
\newcommand{\Gm}{\mathbb{G\ \!\!}_m}
\newcommand{\SB}{\opname{SB}}
\newcommand{\End}{\opname{End}}
\newcommand{\Brauer}{\opname{Br}}
\newcommand{\rep}{\opname{rep}}
\newcommand{\DModule}{\opname{DM}}
\newcommand{\dashto}{\dashrightarrow}
\newcommand{\spec}{\opname{Spec}}
\newcommand{\greatest}{\opname{gr}}
\newcommand{\mCov}{\opname{mCov}}

\input xy
\xyoption{all}

\begin{document}
\title[Multihomogeneous Covariants \& Essential Dimension]{Application of Multihomogeneous Covariants to the Essential Dimension of Finite Groups}
\author{Roland L\"otscher} 
\thanks{The author gratefully acknowledges support from the Swiss National Science Foundation (Schweizerischer Nationalfonds).}
\begin{abstract}
We investigate essential dimension of finite groups over arbitrary fields and give a systematic treatment of multihomogenization, introduced in \cite{KLS}. We generalize the central extension theorem of Buhler and Reichstein, \cite[Theorem 5.3]{BR} and use multihomogenization to substitute and generalize the stack-involved part of the theorem of Karpenko and Merkurjev \cite{KM} about the essential dimension of $p$-groups. One part of this paper is devoted to the study of completely reducible faithful representations. Amongst results concerning faithful representations of minimal dimension there is a computation of the minimal number of irreducible components needed for a faithful representation. 
\end{abstract}

\maketitle

\section{Introduction}
Throughout this paper we work over an arbitrary base field $k$. Sometimes we extend scalars to a larger base field, which will be denoted by $K$. All vector spaces and representations in consideration are finite dimensional over the base field. A quasi-projective variety defined over the base field will be abbreviated as a variety. Unless stated otherwise we will always assume varieties to be irreducible. We denote by $G$ a finite group. A $G$-variety is then a variety with a regular algebraic $G$-action $G \times X \to X, x \mapsto gx$ on it. 
\par
The \emph{essential dimension of $G$} was introduced by Buhler and Reichstein \cite{BR} in terms of \emph{compressions}: A \emph{compression} of a (faithful) $G$-variety $Y$ is a dominant $G$-equivariant rational map $\phi \colon Y\dashto X$, where $X$ is a faithful $G$-variety. \par
\begin{defn}
\label{def:edcomp}
The \emph{essential dimension of $G$} is the minimal dimension of a compression $\phi \colon \AAA(V) \dashto X$ of a faithful representation $V$ of $G$.
\end{defn}

The notion of essential dimension is related to Galois algebras, torsors, generic polynomials, cohomological invariants and other topics, see \cite{BR}. There is a general definition of the essential dimension of a functor from the category of field extensions of $k$ to the category of sets, which is due to Merkurjev, see \cite{Favi}. The essential dimension of $G$ corresponds to the essential dimension of the Galois cohomology functor $K \mapsto H^1(K,G)$. We shall use this only in section \ref{sec:FunctorialEd}. \par
We take the point of view from \cite{KS}, where the \emph{covariant dimension of $G$} was introduced:
A \emph{covariant} of $G$ (over $k$) is a $G$-equivariant ($k$-)rational map $\phi \colon \AAA(V) \dashto \AAA(W)$, where $V$ and $W$ are (linear) representations of $G$ (over $k$). The covariant $\phi$ is called \emph{faithful} if the image of the generic point of $\AAA(V)$ has trivial stabilizer. Equivalently there exists a $\bar{k}$-rational point in the image of $\phi$ with trivial stabilizer.
We denote by $\dim \phi$ the dimension of the closure of the image of $\phi$.
\begin{defn}
\label{def:ed2}
The \emph{essential dimension of $G$}, denoted by $\edim_k G$, is the minimum of $\dim \phi$ where $\phi$ runs over all faithful covariants over $k$. \par
The \emph{covariant dimension of $G$}, denoted by $\cdim_k G$, is the minimum of $\dim \phi$ where $\phi $ runs only over the regular faithful covariants over $k$. \par
\end{defn}
The second definition of essential dimension is in fact equivalent to the first definition, which follows e.g. from \cite[Proposition 2.5]{Fl} or from (the first part of) the following lemma: 
\begin{lem}
\label{le:faithmap}
Let $W$ be a faithful representation of $G$. Then for every affine unirational faithful $G$-variety $X$ there exists a faithful regular $G$-equivariant map $\psi \colon X \to \AAA(W)$. If $X$ contains a $k$-rational point $x_0 \in X(k)$ with trivial stabilizer and $w_0 \in W$ has trivial stabilizer as well, then $\psi$ can be chosen such that $\psi(x_0)=w_0$: 
\end{lem}
\begin{proof}
Choose $f \in k[X]$ such that $f(x_0)=1$ and $f(gx_0)=0$ for $g \neq e$, and define a regular $G$-equivariant map $\psi \colon X \to \AAA(W)$ by
$$\psi(x)=\sum\limits_{g \in G} f(gx) g^{-1} w_0.$$  The map $\psi$ is faithful since $w_0$ is in the image of $\psi$. This shows the second part of the lemma. If $k$ is infinite this immediately implies the first part since in that case the $k$-rational points in $X$ and $\AAA(W)$ are dense. \par
Now let $k$ be a finite field and let $t$ be transcendental over $k$. Since $k(t)$ is infinite we obtain a faithful regular $k(t)$-rational $G$-equivariant map $X_{k(t)} \to \AAA(W\otimes k(t))$ where $X_{k(t)} = X \times_{\spec k} \spec k(t)$ is $X$ with scalars extended to $k(t)$. This corresponds to a homomorphism $W^\ast \otimes k(t) \to k[X]\otimes k(t)$ of representations of $G$ with faithful image, where $W^\ast$ is the dual of $W$ and $k[X]$ is the affine coordinate ring of $X$. Actually we may replace $k[X]\otimes k(t)$ by $U\otimes k(t)$ for some finite-dimensional sub-representation $U\subset k[X]$. By the following Lemma \ref{le:repinj} there exists a homomorphism $W^\ast \to k[X]$ with faithful image, hence a faithful regular $G$-equivariant map $\psi \colon X \to \AAA(W)$. \qed
\end{proof}
\begin{lem}
\label{le:repinj}
Let $W$ and $V$ be (finite-dimensional) representations of $G$ over $k$. Then:
\begin{itemize}
\item If $V \otimes k(t)$ is a quotient of $W \otimes k(t)$ then $W$ is a quotient of $V$. 
\item If $W\otimes k(t)$ injects into $V \otimes k(t)$ then $W$ injects into $V$.
\item If $W \otimes k(t) \to V \otimes k(t)$ is a homomorphism with faithful image, then there exists a homomorphism $W \to V$ with faithful image as well.
\end{itemize}
\end{lem}
\begin{proof}
To show the first claim let $\pi \colon W \otimes k(t) \twoheadrightarrow V\otimes k(t)$ denote the quotient map. Since $t$ is transcendental over $k$ the kernel of $\pi$ can be lifted to a representation $U$ of $G$ over $k$, i.e. $\ker \pi \simeq U\otimes k(t)$. Hence $$(W/U) \otimes k(t) \simeq (W\otimes k(t))/(U\otimes k(t)) \simeq V \otimes k(t).$$ By the theorem of Noether-Deuring this implies $W/U \simeq V$, showing the claim.
The second claim follows from the first claim and dualization. The third claim follows from the first two applied to $V\otimes k(t) \twoheadrightarrow X \otimes k(t)$ and $X\otimes k(t) \hookrightarrow V\otimes k(t)$ where $X$ is a lift of the image of $W\otimes k(t) \to V \otimes k(t)$ to a (faithful) representation of $G$ over $k$. \qed
\end{proof}
We call a faithful regular (resp. rational) covariant \emph{minimal} if $\dim \phi = \cdim_k G$ (resp. $\dim \phi = \edim_k G$). For any faithful representations $V$ and $W$ of $G$ there exists a minimal faithful regular (resp. rational) covariant $\phi \colon \AAA(V) \dashto \AAA(W)$. This is basically another consequence of Lemma \ref{le:faithmap}. At least it shows immediately that the choice of $W$ is arbitrary and if $k$ is infinite one can use $k$-rational points with trivial stabilizer as in \cite[Proposition 2.1]{KS} to show that $V$ can be arbitrarily chosen. For arbitrary fields use e.g.~\cite[Corollary 3.16]{Favi} to see independence of the choice of $V$. \par
\medskip

In sections \ref{sec:mhom} and \ref{sec:properties} we develop the technique of multihomogenization of covariants and derive some of its basic properties. Given $G$-stable gradings $V=\bigoplus_{i=1}^m V_i$ and $W=\bigoplus_{j=1}^n W_j$ a covariant $\phi=(\phi_1,\dotsc,\phi_n) \colon \AAA(V) \dashto \AAA(W)$ is called multihomogeneous if the identities $$\phi_j(v_1,\dotsc,v_{i-1},s v_i,v_{i+1},\dotsc,v_m) = s^{m_{ij}}\phi_j(v_1,\dotsc,v_m)$$ hold true. Here $s$ is an indeterminate and the $m_{ij}$ are integers, forming some matrix $M_{\phi} \in \MMMM_{m\times n}(\ZZZ)$. Thus multihomogeneous covariants generalize homogeneous covariants. A whole matrix of integers takes the role of a single integer, the degree of a homogeneous covariant. It will be shown that the degree matrix $M_{\phi}$ and especially its rank have a deeper meaning with regards to the essential dimension of $G$. Theorem \ref{th:2} states that if each $V_i$ and $W_j$ is irreducible then the rank of the matrix $M$ is bounded from bellow by the rank of a certain central subgroup $Z(G,k)$ (the $k$-center, see Definition \ref{def:kcenter}). Moreover if the rank of $M_{\phi}$ exceeds the rank of $Z(G,k)$ by $\Delta \in \NNN$ then $\edim_k G \leq \dim \phi - \Delta$. This observation shall be useful in proving (partly new) lower bounds to $\edim_k G$ and for most applications in the sequel.
\medskip

In section \ref{sec:exrep} we study faithful representations of $G$, especially faithful representations of small dimension. It is the representation theoretic counterpart to the results on essential dimension obtained in later sections. \par
Section \ref{sec:EsCov} relates essential dimension and covariant dimension. It is well known that the two differ at most by $1$, see the proof of \cite{Re}, which works for arbitrary fields. By generalizing \cite[Theorem 3.1]{KLS} (where $k$ is algebraically closed of characteristic $0$) to arbitrary fields we obtain the precise relation of covariant and essential dimension in case that $G$ has a completely reducible faithful representation. Namely Theorem \ref{th:coves} says that $\cdim_k G = \edim_k G$ if and only if $G$ (is trivial or) has a nontrivial $k$-center, otherwise $\cdim_k G = \edim_k G + 1$.
\par
A generalization of a result from \cite{BR} is obtained in section \ref{sec:CentExt} where the following situation is investigated: $G$ is a (finite) group and $H$ a central cyclic subgroup which intersects the commutator subgroup of $G$ trivially. Buhler and Reichstein deduced the relation $$\edim_k G = \edim_k G/H+1$$ (over a field $k$ of characteristic $0$) for the case that $H$ is a maximal cyclic subgroup of the $k$-center $Z(G,k)$ and has prime order $p$ and that there exists a character of $G$ which is faithful on $H$, see \cite[Theorem 5.3]{BR}. The above theorem was generalized to arbitrary fields in \cite[Theorem 4.5]{Ch}, where for the case of $p = \Chr k>0$ the additional assumption is made that $G$ contains no non-trivial normal $p$-subgroup. Some other partial results were obtained by Brosnan, Reichstein and Vistoli in \cite{BRV} and \cite{BRV2} and by Kraft and Schwarz and the author in \cite{KLS}.
In this paper we give a complete generalization which reads like 
$$\edim_k G = \edim_k G/H + \rank Z(G,k) - \rank Z(G,k)/H,$$ 
where we only assume that $G$ has no non-trivial normal $p$-subgroups if $\Chr k = p>0$ and that $k$ contains a primitive root of unity of high enough order. For details see Theorem \ref{part:centextA}. \par
Section \ref{sec:subpr} contains two additional results about subgroups and direct products, both obtained easily with the use of multihomogeneous covariants.
\par
In section \ref{sec:twist} we shall use multihomogeneous covariants to generalize Florence's twisting construction from \cite{Fl}. The generalized technique gives a substitution for the use of algebraic stacks in the proof of the theorem of Karpenko and Merkurjev about the essential dimension of $p$-groups, which says that the essential dimension of a $p$-group $G$ equals the least dimension of a faithful representation of $G$, provided that the base field contains a primitive $p$-th root of unity. Actually the twisting construction gives more than that. It yields a conjectural formula for the essential dimension of any group $G$ whose socle is central (i.e. such that every nontrivial normal subgroup of $G$ intersects the center of $G$ nontrivially) and whose degrees of irreducible representations satisfy some divisibility property. See Corollary of Conjecture \ref{conjcor:1} for details. \par
In section \ref{sec:FunctorialEd} we consider the situation when multihomogenization fails. This is the case when $G$ does not admit a faithful completely reducible representation. That can only happen if $\Chr k = p>0$ and $G$ contains a nontrivial normal elementary abelian $p$-subgroup $A$. Proposition \ref{pr:estimate} relates the essential dimension of $G$ and $G/A$ by $\edim_k G/A \leq \edim_k G \leq \edim_k G/A + 1$ when $A$ is central.

\section{The technique of multihomogenization}
\label{sec:mhom}
\subsection{Multihomogeneous maps and multihomogenization}
\label{subsec:mhom}
Most of this section can already been found in \cite{KLS}, where multihomogenization has originally been introduced for regular covariants (over $\CCC$). We give a more direct and general approach here. \par
Denote by $X = \Hom(\cdot, \Gm)$ the contravariant functor from the category of commutative algebraic groups (over $k$) to the category of abelian groups, which takes a commutative algebraic group $\Gamma$ to $X(\Gamma)=\Hom(\Gamma,\Gm)$. For example $X(T) = \ZZZ^r$ if $T= \Gm^r$ is a split torus of rank $r=\dim T$. In particular $X(\Gm) = \ZZZ$. \par
Let $T=\Gm^m$ and $T'= \Gm^n$ be split tori. Any homomorphism $D \in \Hom(T,T')$ corresponds to a linear map $X(D) \colon X(T') \to X(T)$ and to a matrix $M_D \in \MMMM_{m \times n}(\ZZZ)$ under the canonical isomorphisms 
$$\Hom(T,T') \cong \Hom(X(T'),X(T)) = \Hom(\ZZZ^n, \ZZZ^m) \cong \MMMM_{m \times n}(\ZZZ)$$
In terms of the matrix $M_D=:(m_{ij})_{1\leq i \leq m,1 \leq j \leq n}$ the homomorphism $D$ is then given by $$D(t_1,\dotsc,t_n)=(t'_1,\dotsc,t'_m) \text{ where } t_j' = \prod_{i=1}^n t_i^{m_{ij}}.$$ \par
The above isomorphisms are compatible with composition of homomorphisms $D \in \Hom(T,T'), D'\in \Hom(T',T'')$ on the left side and multiplication of matrices $M \in \MMMM_{m\times n}(\ZZZ)$, $M' \in \MMMM_{n,r}(\ZZZ)$ on the right side, where $T''$ is another split torus and $r = \rank T''$. That means that $M_{D'\circ D} = M_D \cdot M_{D'}$.

Let $V$ be a vector space equipped with a decomposition $V=\bigoplus_{i=1}^m V_i$. We call $V$ a graded vector space and associate to $V$ the torus $T_V \subseteq \GL(V)$ consisting of those linear automorphisms which are a (non-zero) multiple of the identity on each $V_i$. We identify $T_V$ with $\Gm^m$ acting on $\AAA(V)$ by $$(t_1,\dotsc,t_m)(v_1,\dotsc,v_m) = (t_1v_1,\dotsc,t_m v_m).$$ Let
$W=\bigoplus_{j=1}^n W_j$ be another graded vector space and $T_W\subseteq \GL(W)$ its associated torus.  Let $D \in \Hom(T_V,T_W)$.
A rational map $\phi \colon \AAA(V) \dashto \AAA(W)$ is called $D$-\emph{multihomogeneous} if the diagram
\begin{equation}
\label{eq:multhom}
\xymatrix{
T_V \times \AAA(V) \ar@{->}[rr]^{(t,v)\mapsto tv} \ar@{-->}[d]_{D \times \phi}&& \AAA(V) \ar@{-->}[d]^{\phi} \\
T_W \times \AAA(W) \ar[rr]^{(t',w)\mapsto t'w} && \AAA(W) 
}
\end{equation} commutes. The map $\phi$ is called \emph{multihomogeneous} if it is $D$-multihomogeneous for some $D \in \Hom(T_V,T_W)$.
In terms of the matrix $M_D =: (m_{ij})_{1 \leq i \leq m, 1 \leq j \leq n}$ this means:
\begin{equation}
\label{eq:multhommat}
\phi_j(v_1,\dotsc,s v_i,\dotsc,v_m) = s^{m_{ij}} \phi_j(v_1,\dotsc,v_m), \end{equation} for all $i$ and $j$, as announced in the introduction.
\begin{exa}
\label{ex:1}
Let $V=\bigoplus_{i=1}^m {V_i}$ be a graded vector space. If $h_{ij} \in k(V_i)^\ast$ for $1 \leq i,j \leq m$ are homogeneous rational functions of degree $r_{ij} \in \ZZZ$ then the map
$$ \psi_h \colon \AAA(V) \to \AAA(V),\quad v \mapsto (h_{11}(v_1)\dots h_{m1}(v_m) v_1,\dotsc, h_{1m}(v_1)\dots h_{mm}(v_m) v_m)$$ is multihomogeneous with degree matrix equal to $M_D = (r_{ij} + \delta_{ij})_{1\leq i,j \leq m}$. 
\end{exa}
Let $\phi \colon \AAA(V) \dashto \AAA(W)$ be a multihomogeneous rational map. If the projections $\phi_j$ of $\phi$ to $\AAA(W_j)$ are non-zero for all $j$, then the homomorphism $D \in \Hom(T_V,T_W)$ is uniquely determined by condition \eqref{eq:multhom}. We shall write $D_{\phi}$, $X_\phi$ and $M_{\phi}$ for $D$, $X(D)$ and $M_D$, respectively.
If $\phi_j=0$ for some $j$ then the matrix entries $m_{ij}$ of $M_\phi$ can be chosen arbitrary. Fixing the choice $m_{ij}=0$ for such $j$ makes $M_\phi$ with the property \eqref{eq:multhommat} and the corresponding $D_\phi$ with the property \eqref{eq:multhom} unique again. This convention that we shall use in the sequel has the advantage that adding or removing of some zero-components of the map $\phi$ does not change the rank of the matrix $M_{\phi}$. \par

Given an arbitrary rational map $\phi \colon \AAA(V) \dashto \AAA(W)$ we will produce a multihomogeneous map $H_{\lambda}(\phi)\colon \AAA(V) \dashto \AAA(W)$ which depends only on $\phi$ and the choice of a suitable one-parameter subgroup $\lambda \in \Hom(\Gm,T_V)$. In section \ref{sec:properties} this procedure will be applied to covariants for a group $G$. \par
Let $\nu \colon k(V\times k) = k(s)(V) \to \ZZZ \cup \{\infty\}$ be the discrete valuation belonging to the hyperplane $\AAA(V) \times \{0\} \subset \AAA(V) \times \AAA^1$. So $\nu(0)=\infty$ and for $f \in k(V\times k) \setminus \{0\}$ the value of $\nu(f)$ is the exponent in which the coordinate $s$ appears in a primary decomposition of $f$. Let $O_s \subset k(V\times k)$ denote the valuation ring corresponding to $\nu$. Every $f \in O_s$ can be written as $f=\frac{p}{q}$ with polynomials $p,q$ where $s \nmid q$. For such $f$ we define $\lim{f} \in k(V)$ by $(\lim f)(v) = \frac{p(v,0)}{q(v,0)}$ on the dense open subset in $\AAA(V)$ where $q(v,0) \neq 0$. It is non-zero if and only if $\nu(f)=0$. Moreover $\nu(f - \lim f) > 0$ where $\lim f \in k(V)$ is considered as element of $k(V\times k)$. This follows from writing $(f-\lim(f))(v,s)$
as $$\frac{p(v,s)}{q(v,s)}-\frac{p(v,0)}{q(v,0)} = \frac{p(v,s)q(v,0)-q(v,s)p(v,0)}{q(v,s)q(v,0)},$$
noting that $s$ does not divide the denominator, but $s$ divides the numerator since the numerator vanishes on the hyperplane $\AAA(V) \times \{0\} \subset \AAA(V) \times \AAA^1$. This construction can easily be generalized for rational maps $\psi \colon \AAA(V)\times \AAA^1 \dashto \AAA(W)$ by choosing coordinates on $W$. So for $\psi=(f_1,\dotsc,f_d)$ where $d=\dim W$ and $f_1,\dotsc,f_d \in O_s$ we shall write $\lim \psi$ for the rational map $(\lim f_1,\dotsc,\lim f_d) \colon \AAA(V) \dashto \AAA(W)$. One may check that this definition does not depend on the choice of the basis of $W$.
\par
Let $\lambda \in \Hom(\Gm,T_V)$ be a one-parameter subgroup of $T_V$. Consider $$\tilde{\phi} \colon \AAA(V) \times \Gm \dashto \AAA(W), \quad (v,s) \mapsto \phi(\lambda(s)v)$$ as a rational map on $\AAA(V) \times \AAA^1$.
For $j=1\dots m$ let $\alpha_j$ be the smallest integer $d$ such that all coordinates functions in $s^d \tilde{\phi}_j$ are elements of $O_s$. Actually that works only if $\tilde{\phi}_j \neq 0$. Otherwise we choose $\alpha_j=0$. Let $\lambda' \in \Hom(\Gm,T_W)$ be the one-parameter subgroup corresponding to $\alpha$, i.e. $\lambda'(s) = (s^{\alpha_1},\dotsc, s^{\alpha_n}) \in T_W$ for $s \in \Gm$. Then for $\lambda'(s)\tilde{\phi}(v,s) = \lambda'(s)\phi(\lambda(s)v)$ considered as a rational map $\AAA(V)\times \AAA^1 \dashto \AAA(W)$ we can take its limit: 
$$H_{\lambda}(\phi) = \lim \ \big((v,s) \mapsto \lambda'(s) \phi(\lambda(s)v)\big) \colon \AAA(V) \dashto \AAA(W).$$
The limit $H_{\lambda}(\phi)=(H_{\lambda}(\phi)_1,\dotsc,H_{\lambda}(\phi)_n)$ depends only on $\phi$ and the choice of $\lambda$. By construction we have for $j=1\dots n$: $(H_{\lambda}(\phi))_j \neq 0$ if and only if $\phi_j \neq 0$. \par
It is quite immediate that $H_{\lambda}(\phi)$ is equivariant with respect to the homomorphism of tori $\lambda(\Gm) \to \lambda'(\Gm)$ which sends $\lambda(s)$ to $\lambda'(s^{-1})$. However, to get equivariance for the full tori $T_V$ and $T_W$ we have to choose the one-parameter subgroup $\lambda$ carefully. In any case we have the following

\begin{lem}
For any one-parameter subgroup $\lambda \in \Hom(\Gm,T_V)$ we have $$\dim H_{\lambda}(\phi) \leq \dim \phi.$$
\end{lem}
\begin{proof}
Choose a basis in each $W_j$ and take their union for a basis of $W$. Let $d=\dim W$ and write $\phi=(f_1,\dotsc,f_d)$ with respect to the chosen basis, where $f_j \in k(V)$. Then $H_\lambda(\phi)$ is of the form $(\lim \hat{f}_1,\dotsc, \lim \hat{f}_d)$ where each $\hat{f}_j \in O_s \subset k(V\times k)$ is given by $$\hat{f_j}(v,s) = s^{\gamma_j}f(\lambda(s)v)$$
for some $\gamma_j \in \ZZZ$. 
Choose a maximal subset $S=\{j_1,\dotsc,j_l\}$ of $\{1,\dotsc,d\}$ with the property that $\lim \hat{f}_{j_1},\dotsc,\lim \hat{f}_{j_l}$ are algebraically independent. It suffices to show that $f_{j_1},\dotsc,f_{j_l}$ are then algebraically independent, too. Without loss of generality $j_1=1,\dotsc,j_l=l$. \par
Assume that $f_1,\dotsc,f_l$ are algebraically dependent. Let $p \in k[x_1,\dotsc,x_l] \setminus \{0\}$ with $p(f_1,\dotsc,f_l) = 0$. Since the algebraic independence implies $\lim \hat{f_j} \neq 0$ for $j=1\dots l$ we have $\nu (\hat{f}_j) = 0$.
Set $\gamma=(\gamma_1,\dotsc,\gamma_l)$ and write $p$ in the form 
$$p = \sum_{i \in \ZZZ} p_i \quad \text{where} \quad p_i = \sum_{\beta \in \NNN^l\colon \beta \cdot \gamma = -i} c_\beta x_1^{\beta_1} \dotsb x_l^{\beta_l}.$$
Let $d = \min\{i \in \ZZZ \mid \exists \beta \in \NNN^l \colon \beta \cdot \gamma = - i, c_{\beta} \neq 0\}$. That implies $p_d \neq 0$. For $j=1\dots l$ there exists $\delta_j \in O_s \subset k(V\times k)$ such that $\hat{f_j} - \lim \hat{f_j} = s\delta_j$.
By construction, 
\begin{align*}
0 & =s^{-d}p(f_1,\dotsc,f_l)(\lambda(s) v) = s^{-d}p(s^{-\gamma_1}\hat{f}_1,\dotsc,s^{-\gamma_l}\hat{f}_l)(v) \\
  &= s^{-d}p\big(s^{-\gamma_1}(\lim \hat{f}_1 + s \delta_1),\dotsc,s^{-\gamma_l}(\lim \hat{f}_l + s \delta_l)\big)(v) \\
& = p_d\big(\lim \hat{f}_1,\dotsc,\lim \hat{f}_l\big)(v) + s h(v,s),
\end{align*}
where $h \in O_s$. Taking the limit shows $p_d\big(\lim \hat{f}_1,\dotsc,\lim \hat{f}_l\big)=0,$ which concludes the proof. \qed
\end{proof}

Now the goal is to find a one-parameter subgroup $\lambda \in \Hom(\Gm,T_V)$ such that $H_{\lambda}(\phi)$ becomes multihomogeneous. We can assume that $\phi_j \neq 0$ for all $j$. Write $\phi$ in the form $\phi=\frac{1}{f}(\psi_1,\dotsc,\psi_n)$ where each $\psi_j \colon \AAA(V) \dashto \AAA(W_j)$ is regular and $f \in k[V]$. The space $\Morph(V,W_j)$ of regular maps $\AAA(V) \to \AAA(W_j)$ carries a representation of $T_V$ where $W_j$ is equipped with the trivial action of $T_V$. It decomposes into a direct sum $\Morph(V,W_j) = \bigoplus \Morph(V,W_j)_{\chi}$ taken over all $\chi \in X(T_V)$, where $$\Morph(V,W_j)_{\chi} = \{\psi \in \Morph(V,W_j) \mid \psi(t^{-1}v) = \chi(t) \psi(v) \text{ for all } t \in T_V, v \in \AAA(V)\}.$$ 
Thus $\psi_1,\dotsc, \psi_n$ can be written as a sum $\psi_j = \sum_{\chi} \psi_j^\chi$ where only finitely many $\psi_j^\chi$ are different from $0$. Similarly $f \in k[V] = \Morph(V,k)$ has a decomposition $f=\sum_{\chi} f^\chi$ with the same properties. Let 
$$S(\psi,f) = \{\chi \in X(T_V) \mid f^\chi \neq 0 \text{ or } \exists j: \psi_j^\chi \neq 0\},$$ which is a finite subset of $X(T_V)$.
\begin{lem}
\label{le:ops}
If $T$ is a split torus and $S\subset X(T)$ is a finite subset then there exists a one-parameter subgroup $\lambda \in \Hom(\Gm,T)$ such that the restriction of the map $X(T) \to \Hom(\Gm,\Gm), \chi \mapsto \chi \circ \lambda$ to $S$ is injective.
\end{lem}
\begin{proof}
The claim can easily be shown via induction on the rank $r=\rank T$ of the torus. Identifying $X(T)=\ZZZ^r = \Hom(\Gm,T)$ and $\Hom(\Gm,\Gm)=\ZZZ$ the above map is given by $\ZZZ^r \to \ZZZ, \quad \alpha \mapsto \langle \alpha, \beta \rangle :=\sum_{i=1}^r \alpha_i \beta_i$, where $\beta \in \ZZZ^r$ corresponds to $\lambda$. \qed
\end{proof}
We shall write $\langle \chi, \lambda\rangle$ for the image of $\chi \circ \lambda$ in $\ZZZ$, i.e. $\chi \circ \lambda(s) = s^{\langle \chi , \lambda\rangle}$ for $s \in \Gm$. Now let $\lambda$ be as in Lemma \ref{le:ops} where $T=T_V$ and $S=S(\psi,f)$. Set $\psi_0=f$. Then there are unique characters $\chi_0,\chi_1,\dotsc,\chi_n$ such that $\chi_j \circ \lambda$ is minimal (considered as integer) amongst all $\chi \circ \lambda$ for which $\psi_j^\chi \neq 0$, for each $j=0\dots n$. Then the rational map $\AAA(V)\times \AAA^1\dashto \AAA(W_j)$ (or $\AAA(V)\times \AAA^1 \dashto \AAA^1$ for $j=0$) given by
\begin{align*}
s^{-\langle \chi_j, \lambda\rangle} \psi_j(\lambda(s)v) &= s^{-\langle \chi_j, \lambda\rangle} \sum_{\chi} \psi_j^\chi (\lambda(s)v) \\
&=s^{-\langle \chi_j, \lambda\rangle}  \sum_{\chi} \chi \circ \lambda(s) \psi_j^\chi(v) \\
& = \sum_{\chi} s^{\langle \chi - \chi_j, \lambda\rangle} \psi_j^\chi(v)
\end{align*}
has limit $\psi_j^{\chi_j}$, which implies that $H_{\lambda}(\phi) = \frac{1}{f^{\chi_0}}(\psi_1^{\chi_1},\dotsc,\psi_n^{\chi_n})$.
Define the homomorphism $D \in \Hom(T_V,T_W)$ by
$$D = (\chi_1\chi_0^{-1}, \dotsc,\chi_n\chi_0^{-1}).$$
Then $H_{\lambda}(\phi)(tv) = D(t) H_{\lambda}(\phi)(v)$, showing that $H_{\lambda}(\phi)$ is $D$-multihomogeneous.

\subsection{Existence of minimal multihomogeneous covariants}
\label{subsec:mhomcov}
We now go over to the case where the graded vector spaces $V=\bigoplus_{i=1}^m V_i$ and $W=\bigoplus_{j=1}^n W_j$ are furnished with a representation of $G$. We assume that the tori $T_V$ and $T_W$ commute with the action of $G$ on $V$ and $W$, respectively. Equivalently, the subspaces $V_i$ and $W_j$ are $G$-invariant. We will then represent a covariant $\phi \colon \AAA(V) \dashto \AAA(W)$ as $\phi=\frac{1}{f}\psi$ where $\psi\colon \AAA(V) \to \AAA(W)$ is a regular covariant and $f \in k[V]^G$. For $\lambda \in \Hom(T_V,T_W)$ as in Lemma \ref{le:ops} the rational map $H_{\lambda}(\phi)\colon \AAA(V) \dashto \AAA(W)$ is then multihomogeneous and has dimension $\dim H_{\lambda}(\phi) \leq \dim \phi$. Moreover, $H_\lambda(\phi)$ is again a covariant, since the weight spaces $\Morph(V,W_j)_{\chi}$ and $\Morph(V,k)_{\chi}$ are $G$-stable, so for $j=1\dots n$ the maps $\psi_j^\chi$ and in particular $\psi_j^{\chi_j}$ are covariants for $G$ and the functions $f^\chi$ and in particular $f^{\chi_0}$ are invariants. In general $H_{\lambda}(\phi)$ does not have to be faithful if $\phi$ is. However: 
\begin{lem}
If the representations $W_1,\dotsc,W_n$ are all irreducible, then $H_{\lambda}(\phi)$ is faithful as well.
\end{lem}
\begin{proof}
Let $N_j$ and $N'_j$ denote the stabilizer of the image of the generic point of $\phi_j$ and $H_{\lambda}(\phi_j)$, respectively. It suffices to show $N_j=N'_j$ for $j=1\dots n$. If $\phi_j$ is zero then $H_{\lambda}(\phi_j)=0$ as well and $N_j=G=N'_j$. In the other case both maps are nonzero and their images are $G$-stable subsets of $W_j\otimes k(V)$ spanning $W_j\otimes k(V)$ linearly (since $W_j\otimes k(V)$ is irreducible). Thus $N_j$ and $N'_j$ are both equal to the kernel of the action of $G$ on $W_j$. Again $N_j=N'_j$. \qed
\end{proof}

Thus if we have a minimal faithful covariant $\phi\colon \AAA(V) \dashto \AAA(W)$ and $W=\bigoplus_{j=1}^n W_j$ is a decomposition into irreducible sub-representations, we can always replace it by the multihomogeneous covariant $H_{\lambda}(\phi)$ without loosing faithfulness or minimality. \par

Note that a completely reducible faithful representation $W$ does not exist for every choice of $G$ and $k$. For example if $k=\bar{k}$ and the center of $G$ has an element $g$ of prime-order $p$, then $g$ acts as a primitive $p$-th root of unity on some of the irreducible components of $W$. That is only possible if $\Chr k \neq p$. We use the following:
\begin{defn}
$G$ is called \emph{semi-faithful} (over $k$) if it admits a completely reducible faithful representation (over $k$).
\end{defn}
A criterion for a group to admit a completely reducible faithful representation with any fixed number of irreducible components was given by Shoda \cite{Sh} (in the ordinary case) and Nakayama \cite{Na} (in the modular case). In particular Nakayama obtained \cite[Theorem 1]{Na} that $G$ is semi-faithful over a field of $\Chr k = p>0$ if and only if it has no nontrivial normal $p$-subgroups. One direction follows from Clifford's theorem which says that the restriction of a completely reducible representation to a normal subgroup is again completely reducible and the fact that the only irreducible representation of a $p$-group in characteristic $p$ is the trivial one. For the other implication see Lemma \ref{le:arep}. Therefore we get the following
\begin{cor}
If either $\Chr k = 0$, or $\Chr k = p> 0$ and $G$ has no nontrivial normal $p$-subgroup, there exists a multihomogeneous minimal faithful covariant for $G$.
\end{cor}

\subsection{Multihomogeneous invariants}
Let $V=\bigoplus_{i=1}^m V_i$ be a graded vector space. An element $f \in k(V)$ is called \emph{multihomogeneous} if it is multihomogeneous regarded as a rational map $\AAA(V)\dashto \AAA^1$. Let $G$ be semi-faithful and $V$ a faithful completely reducible representation. The non-zero multihomogeneous invariants form a group under multiplication, denoted by $\MMMM_G(V)$. It is a system of generators for the field $k(V)^G$ of invariants.
\begin{defn}
The \emph{degree module} $\DModule_G(V)$ of $V$ is the submodule of $X(T_V)\simeq\ZZZ^m$ formed by the degrees of multihomogeneous invariants, i.e. the image of the group homomorphism $\deg \colon \MMMM_G(V) \to X(T_V), \quad f \mapsto D_f(\identity_{\Gm})$. Equivalently it is the image of the group homomorphism
$$ \prod_{f \in S} X(\Gm) \to X(T_V)$$ induced by the homomorphisms $X(D_f) \colon X(\Gm) \to X(T_V)$, where $S\subseteq \MMMM_G(V)$ is any finite subset whose degrees generate $\DModule_G(V)$. 
\end{defn}
\begin{defn}
\label{def:kcenter}
The central subgroup $$Z(G,k):=\{g \in Z(G) \mid \zeta_{\ord g} \in k\}$$ of $G$ is called the \emph{$k$-center} of $G$.
\end{defn} 
The $k$-center of $G$ is the largest central subgroup $Z$ for which $k$ contains a primitive root of unity of order $\exp Z$. The groups $Z(G,k)$ and $X(Z(G,k))=\Hom(Z(G,k),\Gm)$ are (non-canonically) isomorphic. The elements of $Z(G,k)$ are precisely the elements of $G$ which act as scalars on every irreducible representation of $G$ over $k$:
\begin{lem}
\label{le:centk}
Let $V=\bigoplus_{i=1}^m V_i$ be any completely reducible faithful representation. Then $\rho_V(Z(G,k)) = T_V \cap \rho_V(G)$.
\end{lem}
\begin{proof}
Since both sides are abelian groups it suffices to prove equality for their Sylow-subgroups. Let $p$ be a prime ($p\neq \Chr k$) and $g \in Z(G)$ be an element of order $p^l$ for some $l\in \NNN_0$. We must show that the following conditions are equivalent:
\begin{enumerate}
\item $g$ acts as a scalar on every $V_i$ 
\item $\zeta_{p^l} \in k$. 
\end{enumerate}
Since $V$ is faithful the order of $g$ equals the order of $\rho(g) \in \GL(V)$, hence the first condition implies the second one. Conversely let $\rho''\colon G \to \GL(V_0)$ be any irreducible representation of $G$. Then the minimal polynomial of $\rho''(g)$ has a root in $k$ since it divides $T^{p^{l}}-1 \in k[T]$ which factors over $k$ assuming the second condition. Hence $\rho''(g)$ is a multiple of the identity on $V'$. In particular this holds for $G \to \GL(V_i)$, proving the claim. \qed
\end{proof}
Degree module and the $k$-center of $G$ are related as follows:
\begin{prop}
\label{pr:coker}
The sequence $\MMMM_G(V) \stackrel{\deg}\longrightarrow X(T_V) \to X(Z(G,k)) \to 1$ is exact and in particular $X(T_V)/\DModule_G(V) \cong X(Z(G,k)) \simeq Z(G,k)$.
\end{prop}
\begin{proof}
Choose a finite subset $S\subseteq \MMMM_G(V)$ such that the degrees of $S$ generate $\DModule_G(V)$.
We may replace the homomorphism $\deg \colon \MMMM_G(V) \to X(T_V)$ by the homomorphism $X(\prod_{f \in S} \Gm) \to X(T_V)$, since they both have image $\DModule_G(V)$. Now the claim becomes equivalent to exactness of the sequence $$1 \to Z(G,k)\to T_V \to \prod_{f \in S} \Gm.$$
Exactness at $Z(G,k)$ follows directly from faithfulness of $V$. Denote by $Q$ the kernel of the last map, which is the intersection of the kernels of the maps $D_f \colon T_V \to \Gm$ taken over all multihomogeneous invariants $f\in S$. Clearly $\rho_V(Z(G,k))\subseteq Q$ because $f$ is $G$-invariant. On the other hand let $\tilde{G}$ be the subgroup of $\GL(V)$ generated by $\rho_V(G)$ and $Q$. Then $\MMMM_G(V) = \MMMM_{\tilde{G}}(V)$ and therefore $k(V)^G = k(V)^{\rho_V(G)} = k(V)^{\tilde{G}}$. This can only happen if $\rho_V(G)=\tilde{G}$. By Lemma \ref{le:centk} this implies $Q=\rho_V(Z(G,k))$, showing the claim. \qed
\end{proof}

Let $\phi=(\phi_1,\dotsc,\phi_n) \colon \AAA(V) \dashto \AAA(W)$ be a faithful multihomogeneous covariant and let $f_1,\dotsc,f_n \in \MMMM_G(V)$ be multihomogeneous invariants. Then $\tilde{\phi}=(f_1\phi_1,\dotsc,f_n \phi_n)\colon \AAA(V) \dashto \AAA(W)$ is again a faithful covariant. That induces an action of the group $\MMMM_G(V)^n$ on the space $\mCov(V,W)$ of multihomogeneous covariants $\AAA(V) \dashto \AAA(W)$, which respects faithfulness. Furthermore we get an action of $\MMMM_G(V)^n$ on the set $S=\{X_\phi \colon \phi \in \mCov(V,W)\} \subseteq \Hom(X(T_W),X(T_V))$ of all degrees associated to multihomogeneous invariants. We will identify the group $\MMMM_G(V)^n$ with the group $\Hom(X(T_W),\MMMM_G(V))$ by associating to an element $\gamma \in \Hom(X(T_W),\MMMM_G(V))$ the $n$-tuple $(f_1,\dotsc,f_n) \in \MMMM_G(V)$ where $f_j = \gamma(\chi_j)$ for the standard basis of $X(T_W)$ formed by the characters $\chi_j\colon T_W \to \Gm, t=(t_1,\dotsc,t_n) \mapsto t_j$. Then the action on degrees is given by 
\begin{align*} 
\Hom(X(T_W),\MMMM_G(V)) \times S & \to S,& \\
(\gamma,s) & \mapsto (\gamma s \colon X(T_W)  \to X(T_V) \\
 & \qquad \chi \mapsto (\deg \gamma(\chi)) \cdot s(\chi).
\end{align*}
From Proposition \ref{pr:coker} we get
\begin{cor}
\label{cor:trans}
The group $\Hom(X(T_W),\MMMM_G(V))$ acts transitively on the set $S$ of all degree matrices associated to multihomogeneous covariants.
\end{cor}
\begin{proof}
Let $s,s' \in S$ and choose $\phi,\phi' \in \mCov(V,W)$ such that $s=X_\phi$ and $s'=X_{\phi'}$. Define $D \in \Hom(T_V,T_W)$ by $D(t)=D_{\phi}(t)D_{\phi'}(t^{-1})$ for $t \in T_V$. Then $D(z)=1$ for all $z \in \rho_V(Z(G,k))$, since $D_{\phi}$ and $D_{\phi'}$ are both the identity on $\rho_V(Z(G,k))$. By Proposition \ref{pr:coker} this is equivalent to saying that $X(D) \in \Hom(X(T_W),\DModule_G(V))$. Therefore $X(D)$ comes from some homomorphism $\gamma \in \Hom(X(T_W), \MMMM_G(V))$. By construction $\gamma s' = s$, finishing the proof. \qed
\end{proof}

Let $\phi \colon \AAA(V) \dashto \AAA(W)$ be a faithful multihomogeneous covariant. Let $N_\phi \in \NNN$ be the greatest common divisor of the entries of the elements of $\image X(D_\phi) \subseteq X(T_V) \cong \ZZZ^m$, where $m=\dim T_V$. Then $N_\phi^{-1} X(D_\phi)\colon X(T_W) \to X(T_V)$ is well defined and its image has a complement in $X(T_V)$. We distinct between two types of elements of $\Hom(X(T_W),\MMMM_G(V))$ relative to $\phi$:
\begin{defn}
A homomorphism $\gamma \colon X(T_W) \to \MMMM_G(V)$ is called of 
\begin{itemize}
\item \emph{type I relative to $\phi$} if it factors through $N_\phi^{-1} X(D_\phi) \colon X(T_W) \to X(T_V)$, i.e. if there exists a commutative diagram of the form
$$\xymatrix{X(T_W) \ar[rr]^{\gamma} \ar[rd]_{N_{\phi}^{-1} X(D_\phi)}& & \MMMM_G(V) \\
& X(T_V) \ar[ru]. & 
}$$
\item \emph{type II relative to $\phi$} if the image of $\gamma$ equals the image of $\ker X(D_\phi) \hookrightarrow X(T_W) \to \MMMM_G(V)$.
\end{itemize}
\end{defn}
\begin{prop}
Every homomorphism $\gamma \colon X(T_W) \to \MMMM_G(V)$ decomposes uniquely as $\gamma=\alpha \cdot \beta$ where $\alpha \colon X(T_W) \to \MMMM_G(V)$ is of type I relative to $\phi$ and $\beta \colon X(T_W) \to \MMMM_G(V)$ is of type II relative to $\phi$.
\end{prop}
\begin{proof}
Uniqueness follows from the fact that the composition $$\ker X(D_\phi) \hookrightarrow X(T_W) \stackrel{N_{\phi}^{-1}X(D_\phi)}\longrightarrow X(T_V)$$ is trivial. It remains to find a decomposition for $\gamma$. Choose decompositions $X(T_W)=\ker X(D_\phi) \oplus A$ and $X(T_V) = \image N_{\phi}^{-1} X(D_\phi) \oplus B$. Define the homomorphisms $\alpha,\beta \colon X(T_W) \to \MMMM_G(V)$ by $$\alpha|_{\ker X(D_\phi)} = 1,\ \beta|_{\ker X(D_\phi)}=\gamma|_{\ker X(D_\phi)} \quad \text{ and } \quad \alpha|_A = \gamma|_A,\ \beta|_A = 1.$$ Clearly $\beta$ is of type II relative to $\phi$ and $\alpha \beta = \gamma$. \par
Note that the homomorphism $N_\phi^{-1} X(D_\phi)\colon X(T_W) \to X(T_V)$ induces an isomorphism from $A$ to its image in $X(T_V)$. Thus we may define $\eps \colon X(T_V) \to \MMMM_G(V)$ by $\eps|_B = 1$ and $\eps(N_{\phi}^{-1} X(D_\phi)(\chi)) = \gamma(\chi)$ for $\chi \in A$. This shows that $\alpha$ is of type I relative to $\phi$, finishing the proof. \qed
\end{proof}
In the sequel the following Lemma will be useful:
\begin{lem}
\label{le:dimcov}
If $\gamma$ is of type I relative to $\phi$ then $\overline{(\gamma\phi)(V_{\bar{k}})}\subseteq \overline{\phi(V_{\bar{k}})}$ and in particular $\dim (\gamma \phi) \leq \dim \phi$. For arbitrary $\gamma$ the dimension of $\gamma\phi$ is at most $\dim \phi + (\rank X(T_W) - \rank M_\phi)$.
\end{lem}
\begin{proof}
Let $\gamma$ be of type I relative to $\phi$. Hence there exists $\eps \colon T_V \to \MMMM_G(V)$ such that $\gamma = \eps \circ N_{\phi}^{-1}X(D_\phi)$. We have rational evaluation maps $\eval_{\gamma} \colon \AAA(V) \dashto T_W$ and $\eval_{\eps} \colon \AAA(V) \dashto T_V$, such that $\eval_{\gamma}(v)=(f_1(v),\dotsc,f_n(v))$ where $f_j$ is the image of the $j$-th standard basis vector under $\gamma$ in $\MMMM_G(V)$, and similiarly for $\eps$. Now let $v \in V_{\bar{k}}$ such that $\eval_{\eps}$ and $\phi$ are defined in $v$. Choose $t \in T_V(\bar{k})$ such that $t^{N_\phi} = \eval_{\eps}(v)$. Then one checks easily that $\eval_\gamma(v) = D_\phi(t)$, whence 
$$(\gamma\phi)(v) = \eval_\gamma(v) \phi(v) = D_\phi(t) \phi(v) = \phi(tv).$$ This proves the first claim.
\par
The second claim follows from the first, since the image of $\ker X(D_\phi) \hookrightarrow X(T_W) \to \MMMM_G(V)$ is generated by $r:=\rank(\ker X(D_\phi)) = \rank X(T_W) - \rank M_\phi$ functions. \qed
\end{proof}

\section{Properties of multihomogeneous covariants}
\label{sec:properties}
\subsection{The rank of the degree-matrix of a multihomogeneous covariant}
Let $G$ be semi-faithful and $V=\bigoplus_{i=1}^m V_i, W=\bigoplus_{j=1}^n W_j$ be two faithful representations of $G$.
For a faithful multihomogeneous covariant $\phi \colon \AAA(V) \dashto \AAA(W)$ we will prove the following interpretation of the rank of the degree-matrix $M_\phi$: 

\begin{theorem}
\label{th:2}
Let $\phi\colon \AAA(V) \dashto \AAA(W)$ be a faithful multihomogeneous covariant. Assume that $W_1,\dotsc,W_n$ are irreducible. 
$$\edim_k G - \rank Z(G,k) \leq \dim \phi - \rank M_{\phi}.$$ 
If furthermore $V_1,\dotsc,V_n$ are irreducible then 
$$\rank M_{\phi} \geq \rank Z(G,k)$$ 
with equality if $\phi$ is minimal.
\end{theorem}

\begin{proof}
Let $Z:=Z(G,k)$. We first prove the second inequality. Since $\phi$ is at the same time equivariant with respect to the tori- and $G$-action $g\phi(v) = \phi(gv) = (D_\phi g)\phi(v)$ for $g \in Z$. Thus the map $D_\phi$ is the identity restricted to $Z$. This implies $Z = D_\phi(Z) \subset D_\phi(T)$, whence $\rank M_\phi = \dim D_\phi(T) \geq \rank Z$.
The first inequality follows from the following: \qed
\end{proof}
\begin{prop}
\label{pr:4}
Let $\phi = (\phi_1,\dotsc,\phi_n) \colon \AAA(V) \dashto \AAA(W)$ be a faithful rational multihomogeneous covariant. Assume that each $W_j$ in the decomposition of $W$ is irreducible. If $\rank M_{\phi} \geq \rank Z(G,k)$ there exists a sub-torus $S\subseteq D_\phi(T_V)$ of dimension $\rank M_{\phi} - \rank Z(G,k)$ and a $G$-invariant open subset $W' \subseteq \AAA(W)$ on which $D_\phi(T_V)$ acts freely such that the action of $G$ on the quotient $(\overline{\image \phi} \cap W')/S$ is faithful.
\end{prop}
\begin{proof}
Let $Z:=Z(G,k)$. The torus $D_\phi(T_V)$ has dimension $d:=\rank M_{\phi} \geq r:=\rank Z$. By the elementary divisor theorem there exist integers $c_1,\dotsc,c_r > 1$ and a basis $\chi_1,\dotsc,\chi_d$ of $X(D_\phi(T_V))$ such that $$Z = \bigcap_{i=1}^r \ker \chi_i^{c_i} \cap \bigcap_{j=r+1}^d \ker \chi_j.$$ 
Set $S:= \bigcap_{i=1}^r \ker \chi_i$. This is a subtorus of $D_\phi(T_V)$ of rank $d-r= \rank M_{\phi} - \rank Z$ with $S \cap Z = \{1\}$.
\par Let $W':=\prod_{j=1}^n W'_j$, where $W'_j := \AAA(W_j) \setminus \{0\}$ if $\phi_j \neq 0$ and $W'_j:=\AAA(W_j)$ otherwise. Our convention that $(M_\phi)_{ij}=0$ if $\phi_j=0$ implies that $D_{\phi}(T_V)$ (and therewith $S$) acts freely on $W'$. Let $X:=\overline{\image(\phi)}$ and set $X':=X\cap W'$. Let $\pi \colon \AAA(W) \dashto W'/S$ be the rational projection map. The kernel of the action of $G$ on $X'/S$ is contained in $Z(G,k)$ by the next lemma. Since $Z(G,k)\cap S = \{e\}$ it is trivial.
Hence the rational map $\AAA(V) \stackrel{\phi}\dashrightarrow X \dashto X'/S$ is a compression and $\edim_k G \leq \dim X'/S = \dim X - \dim S = \dim \phi - (\rank M_{\phi} - \rank Z)$. \qed
\end{proof}
\begin{lem}
\label{le:keract}
Let $\phi \colon \AAA(V) \dashto \AAA(W)$ be a faithful multihomogeneous covariant. Let $P:=\prod\limits_{j:\phi_j \neq 0} \PPP(W_j) \ \times \prod\limits_{j:\phi_j = 0} \AAA(W_j)$ and $\pi \colon \AAA(W) \dashto P$ the obvious $G$-equivariant rational map and let $X:=\overline{\image \phi}$. Then the kernel $Q$ of the action of $G$ on $\pi(X)$ equals $Z(G,k)$. 
\end{lem}
\begin{proof}
The elements of the $k$-center $Z(G,k)$ act as scalar on $\AAA(W_j)$ for each $j$. This implies that $Z(G,k)$ is contained in $Q$. Conversely let $g \in Q$ and fix some $j \in \{1,\dotsc,n\}$ with $\phi_j\neq 0$. We want to show that $g$ acts by multiplication of a (fixed) scalar on $W_j$. From this the inclusion $Q \subseteq Z(G,k)$ follows, since $\bigoplus_{j:\phi_j \neq 0} W_j$ is already a faithful completely reducible representation of $G$. \par
Let $Y\subseteq \AAA(W_j)$ denote the projection of $X \cap W'$ to $\AAA(W_j)$. Since $g$ acts trivially on $\pi(X)$ there exists for every field extension $k'/k$ and $y \in Y(k')$ some $\lambda_y \in \Gm(k')$ such that $gy=\lambda_y y$. Since $g$ has only finitely many eigenvalues $\alpha_1,\dotsc,\alpha_r \in \Gm(\bar{k})$ the same holds for its closure $\overline{Y}$. Moreover since $\overline{Y}$ is irreducible the scalar $\lambda:=\lambda_y$ does not depend on $y$. Since $\phi_j \neq 0$ the variety $\overline{Y}$ contains a non-zero $k(V)$-rational point $y_0$. By irreducibility of $W_j\otimes k(V)$ the set $\{g'y_0\mid g' \in G\}$ spans $W_j \otimes k(V)$ as a $k(V)$-vector space. It follows that $g$ acts by multiplication of $\lambda$ on $W_j \otimes k(V)$, hence in the same manner on $W_j$, which completes the proof. \qed
\end{proof}

To illustrate the usefulness of the existence of minimal faithful multihomogeneous covariants and Lemma \ref{le:keract} we give a simple corollary. Its first part was already established in \cite[Theorem 6.1]{BR}.
\begin{cor}
\label{cor:ab}
Let $A$ be abelian and assume that $k$ contains a primitive root of unity of order $\exp A$. Then $$\edim_k A = \rank A.$$ If $G$ is semi-faithful and if $\edim_k G \leq \rank Z(G,k) + 1$, then $G$ is an extension of a subgroup of $\PGL_2(k)$ by $Z(G,k)$. \par
If $\edim_k G \leq \rank Z(G,k)$ then $G=Z(G,k)$, hence abelian with $\zeta_{\exp G} \in k$.
\end{cor}
\begin{proof}
The inequality $\edim_k A \leq \rank A$ is easy to see, because $A$ has a faithful representation of dimension $\rank A$. Let $V$ be a completely reducible faithful representation of $G$ and let $\phi \colon \AAA(V)\dashto \AAA(V)$ be a minimal faithful multihomogeneous covariant of $G$. We may assume that $\phi_j \neq 0$ for all $j$. The group $G/Z(G,k)$ then acts faithfully on the image of $\pi_V \circ \phi \colon \AAA(V) \dashto \AAA(V) \dashto \PPP(V)$, which has dimension at most $\dim \phi - \rank Z(G,k) = \edim_k G - \rank Z(G,k) \leq 1$. Thus $G/Z(G,k)$ embeds into $\PGL_2(k)$. Now if $\edim_k G \leq \rank Z(G,k)$ then $\edim_k G = \rank Z(G,k)$ and the image of $\pi_V \circ \phi$ must be a point, whence $G=Z(G,k)$. \qed
\end{proof}
\begin{rem}
The second part of Corollary \ref{cor:ab} can be used to classify semi-faithful groups with $\edim_k G - \rank Z(G,k) \leq 1$. For example if $\edim_k G \leq 2$ and $Z(G,k)$ is nontrivial one should obtain with the arguments of \cite[section 10]{KS} that $G \hookrightarrow \GL_2(k)$. We haven't checked that in detail, but one observes that the additional possibilities for subgroups of $\PGL_2(k)$ arising in positive characteristic are not semi-faithful.
\end{rem}

\subsection{Behavior under refinement of the decomposition}
Let $V=\bigoplus_{i=1}^m V_i$ be a graded vector space. For each $i$ let $V_i=\bigoplus_{k=1}^{d_i} V_{ik}$ be a grading of $V_i$. We call the grading $V=\bigoplus_{i,k} V_{ik}$ a \emph{refinement} of the grading $V=\bigoplus_i V_i$. Let $\phi=(\phi_1,\dotsc,\phi_n) \colon \AAA(V) \dashto \AAA(W)$ be a multihomogeneous rational map. We consider refinements both in $V$ and in $W=\bigoplus_{j=1}^n W_j$ where $W_j = \bigoplus_{l=1}^{e_j} W_{jl}$, and study the behavior of the rank of the degree matrix. Set $d=\sum_{i=1}^m d_i$ and $e=\sum_{j=1}^n e_j$.
\begin{prop}
\label{pr:refin}
\begin{enumerate}
\item Refinement in $V$: Let $\lambda$ be a one-parameter subgroup of $T_V = \Gm^d$ such that $H_{\lambda}(\phi) \colon \AAA(V) \dashto \AAA(W)$ is multihomogeneous w.r.t. the refined grading on $V$ and the old grading on $W$. Then $$\rank M_{H_{\lambda}(\phi)} \geq \rank M_\phi.$$
\item Refinement in $W$: The map $\phi$ can be considered as a multihomogeneous map $\phi'\colon \AAA(V)\to \AAA(W)$ with respect to the gradings $V=\bigoplus_i V_i$ and $W=\bigoplus_{jl} W_{jl}$, where
$$\rank M_{\phi'} = \rank M_\phi.$$
\item Refinement in both $V$ and $W$: Consider $\phi'$ as above and let $\lambda$ be a one-parameter subgroup of $T_V = \Gm^d$ be such that $H_\lambda(\phi') \colon \AAA(V) \dashto \AAA(W)$ is multihomogeneous w.r.t. the refined grading on both $V$ and $W$. Then 
$$\rank M_{H_\lambda(\phi')} \geq \rank M_\phi.$$
\end{enumerate}
\end{prop}
\begin{proof}
\begin{enumerate}
\item Let $(a_{i,j}) = M_\phi \in \MMMM_{m,n}(\ZZZ)$ and $(b_{ik,j}) = M_{H_{\lambda}(\phi)} \in \MMMM_{d,n}(\ZZZ)$ be the degree matrices of $\phi$ and $H_{\lambda}(\phi)$, respectively. Since $H_{\lambda}(\phi)$ is still multihomogeneous with respect to the old decomposition of $V$ we have $\sum_{k=1}^{d_i} b_{ik,j} = a_{i,j}$ for $i=1\dots m$ and $j=1\dots n$. Therefore the span of the rows of $M_{H_{\lambda}(\phi)}$ contains the span of the rows of $M_{\phi}$. Hence $\rank M_{H_{\lambda}(\phi)} \geq \rank M_{\phi}$.
\item The maps $\phi_{jl} \colon V \dashto W_{jl}$ are still multihomogeneous of the same degree as $\phi_j \colon \AAA(V) \dashto \AAA(W_j)$, as long as they are non-zero. If $\phi_j$ is non-zero then also one of the $\phi_{jl}$ for $l=1\ldots e_j$. Recall that by convention the matrix entries for zero-components are zero, so that they do not influence the column span of the matrix. Thus the column span of $M_\phi$ equals the column span of $M_{\phi'}$ and hence $\rank M_\phi = \rank M_{\phi'}$.
\item follows from $(A)$ and $(B)$.
\end{enumerate} \qed
\end{proof}

\section{Completely reducible faithful representations}
\label{sec:exrep}
\subsection{Minimal number of irreducible components}
\label{subsec:mincomp}
In this section we will compute the minimal number of irreducible components of a faithful representation of any semi-faithful group. As a consequence we obtain a characterization of groups, which have a faithful representation with any fixed number of irreducible components. Groups admitting an irreducible faithful representation over an algebraically closed field of characteristic $0$ have been characterized in \cite{Ga}. A criterion for a group to admit a faithful representation with any fixed number of irreducible components was given by Shoda \cite{Sh} (in the ordinary case) and Nakayama \cite{Na} (in the modular case). Their criterion is formulated in a way quite different from Gaschuetz's and our characterization. \par

\begin{defn}
A \emph{foot} of $G$ is a minimal nontrivial normal subgroup of $G$.
The subgroup of $G$ generated by the (abelian) feet of $G$ is called the (abelian) \emph{socle of $G$}, denoted by $\soc(G)$ (resp. $\socAb(G)$).
\end{defn}
By construction $\soc(G)$ and $\socAb(G)$ are normal. The following Lemma is well known  and a generalization to countable groups can be found in \cite{BH}.
\begin{lem} 
\label{le:soc}
$\soc(G) = \socAb(G) \times N_1 \times \dotsb \times N_r$, where $N_1,\dotsc,N_r$ are all the non-abelian feet of $G$.
\end{lem}
For a $\ZZZ G$-module $A$ denote by $\grk(A)$ the minimum number of generators:
$$\grk(A):=\min\left\{r \in \NNN_0 \mid \exists a_1,\dotsc,a_r \in A: \left\langle a_1,\dotsc,a_r \right\rangle_{\ZZZ\! G} = A \right\} \in \NNN_0.$$
\begin{prop}
\label{th:mcrep}
Let $G$ be a semi-faithful group. Then the minimal number of factors of a decomposition series of a faithful representation of $G$ over $k$ equals $\grk \socAb(G)$ if $\socAb(G) \neq \{e\}$ and $1$ if $\socAb(G)$ is trivial. Moreover the minimum is attained by a completely reducible representation.
\end{prop}
We start with a lemma explaining how to pass from arbitrary to completely reducible representations.
\begin{lem}
\label{le:arep}
Let $V$ be a faithful representation of $G$ and $\FFFF=V=V_1 \supsetneq V_2 \supsetneq \dotsb \supsetneq V_r \supsetneq V_{r+1}=\{0\}$ be a $G$-stable flag. Assume that $\Chr k = p>0$. If $G$ does not contain a nontrivial normal subgroup of $p$-power order then the associated graded representation $\greatest_{\FFFF} V = \bigoplus_{i=1}^r V_{i}/V_{i+1}$ is faithful as well. In particular such a group $G$ is semi-faithful (over $k$).
\end{lem}
\begin{proof}
It is well known that an element of finite order in a unipotent group in characteristic $p$ has $p$-power order. Therefore the kernel of the representation $\greatest_{\FFFF} V$ is a normal subgroup of $G$ of $p$-power order, which by assumption must be trivial. The last statement follows from taking for $\FFFF$ a decomposition series. \qed
\end{proof}
For the proof of Proposition \ref{th:mcrep} we work with two lattices: Set $A:=\socAb(G)$ and let $A^\ast:=\Hom(A,\bar{k}^\ast)$ denote its group of characters over $\bar{k}$, which is again a $\ZZZ G$-module by endowing $\bar{k}^\ast$ with the trivial $G$-action. Denote by $L(A)$ and $L(A^\ast)$ the lattices of $\ZZZ G$-invariant subspaces of $A$ and $A^\ast$, respectively, where the meet-operation is given by $B \cap C$ and the join-operation by $B \cdot C$.

\begin{lem}
\label{le:latiso}
Assume that either $\Chr k = 0$ or $\Chr k = p > 0$ and $p \nmid |A|$.
\begin{enumerate}
\item The map $$\alpha \colon L(A^\ast) \to L(A), \quad \LLL \mapsto \{a \in A \mid \ell(a) = 1 \ \forall \ell \in \LLL\}$$ yields an anti-isomorphism of $L(A^\ast)$ and $L(A)$ with inverse given by $\alpha^{-1}(B)= \{\ell \in A^\ast \mid \ell(a) =1 \ \forall a \in A\}$.
\item There exists a (non-canonical) isomorphism of lattices $$\beta \colon L(A) \stackrel{\simeq}{\longrightarrow} L(A^\ast)$$ which preserves size, i.e. $|\beta(B)| = |B|$ for all $B \in L(A)$.
\item $\grk(A)=\grk(A^\ast)$.
\end{enumerate}
\end{lem}
\begin{proof}
 \begin{enumerate}
\item The proof is straightforward. 
\item The $\ZZZ G$-module $A$ is semi-simple by construction and thus decomposes into isotypic components. Every submodule of $A$ is isomorphic to the direct sum of its intersections with the isotypic components and it suffices to show the claim for every isotypic component of $A$. Thus assume $A = (\FFF_q)^m \otimes V$ for some prime $q \neq \Chr k$, some natural number $m$ and some irreducible $\FFF_q G$-module $V$, where $(\FFF_q)^m$ is equipped with the trivial action of $G$. Hence we may identify $A^\ast = (\FFF_q)^m \otimes V^\ast$. Every $\ZZZ G$-invariant subgroup of $A$ is now of the form $W \otimes V$ for some sub vector-space $W \subset \FFF_q^m$. Define $\beta \colon L(A) \to L(A^\ast)$ by $\beta (W \otimes V) = W \otimes V^\ast$. Then $\beta$ is an isomorphism of lattices and preserves size, since the assumption $p \nmid |A|$ implies $|V^\ast|=|V|$.
\item Let $E_r \subseteq A$ for $r \in \NNN$ denote the (possibly empty) set of generating $r$-tuples of the $\ZZZ G$-module $A$ and let $\max(L(A))$ be the set of maximal non-trivial elements of $L(A)$. The two sets are related by:
$$E_r = A^r \setminus \bigcup\limits_{M \in \max(L(A))} M^r.$$
Similarly for $E^\ast_r \subseteq A^\ast$ and $\max(L(A^\ast))$ defined correspondingly with $A^\ast$ in place of $A$ we have 
\begin{align*}
E^\ast_r & = (A^\ast)^r \setminus \bigcup\limits_{\LLL \in \max(L(A^\ast))} \LLL^r \\
	& = (\beta(A))^r \setminus \bigcup\limits_{M \in \max(L(A))} (\beta(M))^r 
\end{align*}
We claim for any $r$ that $|E_r|=|E^\ast_r|$. This implies in particular that $A$ is generated by $r$ elements if and only if $A^\ast$ is, hence $\grk(A) = \grk(A^\ast)$. The claim follows from part (B) and the exclusion principle, which says that for subsets $Y_1,\dotsc,Y_t$ of a set $Y$ we have 
$$|Y \setminus \cup_{i=1}^t Y_i| = |Y|-\sum_{i=1}^t (-1)^{t+1} \sum_{\nu_1< \dotsb <\nu_i } |Y_{\nu_1} \cap \dotsb  \cap Y_{\nu_i}|
$$
\end{enumerate} \qed
\end{proof}
For the case that $k$ is not algebraically closed, we need to deal with irreducible representations which are not absolutely irreducible:
\begin{lem}
\label{le:elab}
\begin{enumerate}
\item Let $q \neq \Chr k$ be a prime and $A$ be an elementary abelian $q$-group. Then each non-trivial irreducible representation of $A$ (over $k$) is isomorphic to a sub-representation of
$$V_{\langle \chi \rangle }:=\left\{\sum\limits_{C \in A/\ker \chi} \gamma_C (\sum\limits_{a \in C} a) \in kA \mid \sum\limits_{C \in A/ker \chi} {\gamma_C} = 0 \right\}$$
where $\chi \in \Hom(A,\bar{k}^\ast)$, $\chi \neq 1$.
\item Let $A =\bigoplus_{i=1}^m A_{q_i}$ where $q_1,\dotsc,q_m \neq \Chr k$ are distinct primes and $A_{q_i}$ is an elementary abelian $q_i$ group. Then every irreducible representation $V$ of $A$ is an exterior tensor product of irreducible representations of $A_{q_1}, \dotsc, A_{q_m}$. Let $\chi_1,\dotsc,\chi_r$ be the characters appearing in $V \otimes_k \bar{k}$. Then $\langle \chi_1,\dotsc,\chi_r \rangle = \langle \chi_i\rangle $ for every $i=1\dots r$.
\end{enumerate}
\end{lem}
\begin{proof}
\begin{enumerate}
\item It suffices to show that the group algebra $kA$ decomposes as $$\bigoplus\limits_{\langle \chi \rangle \subseteq \Hom(A,\bar{k}^\ast)} V_{\langle \chi \rangle},$$ where we set $V_{\langle \chi \rangle} = k \sum_{a \in A} a$ for $\chi=1$, which has dimension one. Let $n:=\dim_{\FFF_q} A$. There are precisely $\frac{q^n-1}{q-1}$ nontrivial subgroups of the form $\langle \chi \rangle$ and the corresponding subspaces $V_{\langle \chi \rangle}$ all have dimension $q-1$. Since $(q-1)\cdot \frac{q^n-1}{q-1} + 1 \cdot 1 = q^n = |A| = \dim_k kA$ it remains to show that the subspaces $V_{\langle \chi \rangle }$ form a direct sum, for which we may pass to an algebraic closure. Consider the elements $\eps_{\chi}:=\sum_{a \in A} \chi(a^{-1}) a \in \bar{k}A$ for $\chi \in \Hom(A,\bar{k}^{\ast})$, which are $\bar{k}$-linearly independent. Then $V_{\langle \chi \rangle} \otimes \bar{k}$ has $\bar{k}$-basis $\eps_{\chi},\dotsc,\eps_{\chi^{q-1}}$ for $\chi \neq 1$ and $V_{\langle 1 \rangle}$ has basis $\eps_0$. That shows the claim.
\item Writing $kA = kA_{q_1} \otimes \dotsb \otimes kA_{q_m}$ the first claim follows from the fact that the group algebras $kA_{q_i}$ are of coprime dimensions. The second claim follows now from the description in (A), noting that the representation $V_{\langle \chi \rangle}$ has character $\sum_{i=1}^{q-1} \chi^i$.
\end{enumerate} \qed
\end{proof}
The following lemma contains the crucial observation for our study of faithful representations.
\begin{lem}
\label{le:socrep}
Let $V=\bigoplus_{i=1}^m V_i$ be a representation of $G$ with each $V_i$ irreducible. Let $A:=\socAb(G)$ and choose for every $i$ some character $\chi_i \in A^\ast$ appearing in $V_i|_A \otimes \bar{k}$. Then $V$ is faithful if and only if the characters $\chi_1,\dotsc,\chi_m$ generate $A^\ast$ as a $\ZZZ G$-module and no nonabelian foot of $G$ is in the kernel of $V$.
\end{lem}
\begin{proof}
Let $\LLL:= \langle \chi_1, \dotsb, \chi_m \rangle_{\ZZZ G} \in L(A^\ast)$. 
Assume that $\LLL \neq A^\ast$. Let $\alpha$ be the lattice anti-isomorphism from Lemma \ref{le:latiso}(A) and set $B:=\alpha(\LLL) \subseteq A$, which is then a non-trivial normal subgroup of $A$ contained in the kernel of each $\chi_i$ and of any power of $\chi_i$. Let $W_i$ be any irreducible sub-representation of $V_i|_A$ containing the character $\chi_i$. By Lemma \ref{le:elab} $W_i \otimes \bar{k} = \sum \bar{k}_{\chi_i^{\alpha_{ij}}}$ for some $\alpha_{ij} \in \NNN$. Therefore $B$ acts trivially on $W_i$. Now since $V_i$ is irreducible, $V_i = \sum_{g \in G} g W_i$ as vector spaces. For $b \in B$ and $w \in W_i$ we have $bgw = g(g^{-1}bg)w = gw$, since $B$ is normal. Thus $B$ acts trivially on $V$. Hence $V$ is not faithful. \par
Conversely assume that $V$ is not faithful and no noabelian foot of $G$ is in the kernel of $V$. Hence some abelian foot $B$ is in the kernel of $V$. This implies that $B$ lies in the kernel of each $\chi_i$, whence in the kernel of each element of $\LLL$. This implies that $\LLL \neq A^\ast$. \qed
\end{proof}

Now we are ready for the proof of the proposition.
\begin{proof}[Proof of Proposition \ref{th:mcrep}]:
Recall that a group admitting a nontrivial normal subgroup of $p$-power order is not semi-faithful in characteristic $p$. From now on assume that $p \nmid |A|$ where $A:=\socAb(G)$. \medskip
\begin{itemize}
\item["$\geq$"] Let $V$ be a faithful representation of $G$ over $k$. We want to show that the number of factors of a decomposition series of $V$ is at least the maximum of $\grk(A)$ and $1$. Clearly it is at least $1$. By Lemma \ref{le:arep} we may assume that $V$ is completely reducible. Lemma \ref{le:socrep} implies that the number of irreducible components of $V$ is at least $\grk(A^\ast)$, which equals $\grk(A)$ by Lemma \ref{le:latiso}(C).
\medskip
\item["$\leq$"] We must construct a faithful representation $V$ over $k$ with at most $\grk(A)$ irreducible components if $A$ is non-trivial, and a faithful irreducible representation $V$ over $k$ if $A$ is trivial. We first reduce to the case of $k$ being algebraically closed: Assume that $\bigoplus_{i=1}^n V_i$ is a decomposition of a faithful representation into irreducible representations over $\bar{k}$. For each $i$ take any irreducible representation $V_i'$ over $k$ which contains $V_i$ as a decomposition factor over $\bar{k}$. Then $\bigoplus_{i=1}^n V_i'$ is a faithful representation over $\bar{k}$ and has the same number of irreducible components. \par
Let $N_1,\dotsc,N_t$ be the non-abelian feet of $G$. By Lemma \ref{le:soc} the socle of $G$ decomposes as $\soc G=A\times N_1\times \dotsc \times N_t$. For each $i$, since $N_i$ has composite order it has a nontrivial irreducible representation $W_i$. The (exterior) tensor product $W:=W_1 \otimes \dotsb \otimes W_t$ is then irreducible (since $k=\bar{k}$) and does not contain any of $N_1,\dotsb,N_t$ in its kernel. If $A$ is trivial this gives an irreducible representation of $\soc G$ with the property that no foot of $G$ is contained in its kernel. Any irreducible representation whose restriction to $\soc G$ contains $W$ is then faithful. \par
From now on assume $A$ to be non-trivial. There exist $r:=\grk(A^\ast)= \grk(A)$ characters $\chi_1,\dotsc,\chi_r$ of $A$ which generate the $\ZZZ G$-module $A^\ast$. For every $i$ choose an irreducible representation $V_i$ of $G$ whose restriction to $\soc G$ contains the irreducible representation $k_{\chi_i} \otimes W$. Set $V:=\bigoplus_{i=1}^r V_i$. By Lemma \ref{le:socrep} the representation $V$ is faithful. Moreover it has the required number of irreducible components. This finishes the proof.
\end{itemize} \qed
\end{proof}
\begin{rem}
\label{rem:asgr}
The situation for non-semi-faithful groups is completely different, in so far that the abelian socle tells nothing about the number of decomposition factors needed for a faithful representation. Take for example the groups $\ZZZ/p^n\ZZZ$, $n\geq 1$, whose abelian socle are all isomorphic although for large $n$ these groups need more than any fixed number of decomposition factors for a faithful representation.
\end{rem}
\begin{rem}
More generally let $\Gamma$ be any subgroup of $\Aut(G)$ containing the inner automorphisms. One can define $\Gamma$-faithful representations, $\Gamma$-feet, $\Gamma$-socle, abelian $\Gamma$-socle (denoted in the sequel by $A^\Gamma(G)$) as in \cite{BH} and generalize Proposition \ref{th:mcrep} in the following way: If $\Chr k = 0$ or $\Chr k = p>0$ and $p \nmid |A^\Gamma(G)|$ then the minimal number of irreducible components of a completely reducible $\Gamma$-faithful representation of $G$ equals the maximum of $\rank_{\ZZZ \! \Gamma} A^\Gamma(G)$ and $1$. The proof remains basically the same.
\end{rem}

There is the following application:
\begin{cor}
Let $n \in \NNN$ and $H\subseteq G$ be a subgroup containing $\socAb(G)$ and assume that $H$ has a faithful representation over $k$ with $n$ decomposition factors. If $\Chr k \nmid |\socAb(G)|$ then $G$ has a faithful representation with $n$ decomposition factors as well.
\end{cor}
\begin{proof}
This is a consequence of the following Lemma \ref{le:socs} together with Proposition \ref{th:mcrep}. Observe that $\Chr k \nmid  |\socAb(G)|$ implies that $\Chr k \nmid |\socAb(H)|$, hence both groups are semi-faithful. \qed
\end{proof}
\begin{lem}
\label{le:socs}
If $H\subseteq G$ is a subgroup containing $\socAb(G)$ then $\rank_{\ZZZ \! H}\socAb(H) \geq \grk \socAb(G)$.
\end{lem}
\begin{proof}
Let $h_1,\dotsc,h_r$ generate $\socAb(H)$ as a $\ZZZ H$-module, where $r=\rank_{\ZZZ \! H}(\socAb(H))$. Let $N$ be an $H$-invariant complement of $\socAb(H) \cap \socAb(G)$ in $\socAb(H)$. Write $h_i = (g_i,n_i)$ where $n_i \in N$ and $g_i \in \socAb(H) \cap \socAb(G)$. Then $g_1,\dotsc, g_r$ generate $\socAb(H) \cap \socAb(G)$ as a $\ZZZ H$-module. We show that $g_1,\dotsc,g_r$ generate $\socAb(G)$ as a $\ZZZ G$-module, which gives the claim. Let $A$ be any abelian foot of $G$. By assumption $A\subseteq \socAb(G) \subseteq H$. Let $B \subseteq A$ be a $H$-foot.
 By construction $B \subseteq \socAb(H) \cap \socAb(G)$, which is generated by $g_1,\dotsc,g_r$ as a $\ZZZ H$-module. Since $A$ is minimal, the $\ZZZ G$-module generated by $B$ equals $A$. Hence $A$ is contained in the $\ZZZ G$-module generated by $g_1,\dotsc,g_r$. Since this holds for every abelian foot $A$ of $G$ the claim follows. \qed
\end{proof}
There is a simple lower bound on the number of irreducible components needed for a faithful representation, namely the rank of the center of $G$. Since representations for which the bound is reached are of some special interest later, we give it a name:
\begin{defn}
A faithful representation $V$ of a semi-faithful group $G$ is called \emph{saturated} if it is the direct sum of $\rank Z(G)$ many irreducible representations of $G$. \par
The group $G$ is called \emph{saturated} if it has a (faithful) saturated representation. Equivalently (by Proposition \ref{th:mcrep}): $$\rank Z(G) = \grk \socAb(G) \geq 1.$$
\end{defn}
It is sometimes advantageous to pass to saturated groups by taking the product with cyclic groups of high enough rank:
\begin{prop}
\label{pr:satrep}
Let $\ell\neq \Chr k$ be any prime number such that $\zeta_\ell \in k$. Assume that $G$ has a completely reducible faithful representation $V=\bigoplus_{i=1}^n V_i$, each $V_i$ irreducible. Let $r$ be the rank of the $\ell$-Sylow subgroup of $Z(G)$. Then $V$ carries a faithful representation of $G\times C_{\ell}^{n-r}$.
\end{prop}
\begin{proof}
We proceed by induction on $n-r$. If $n-r=0$ there is nothing to show. Otherwise $r < n$ and there exists $i\in \{1,\dotsc,n\}$ such that no element of $G$ acts by multiplication of a primitive $\ell$-th root of unity on $V_i$ and trivially at the same time on every $V_j$ for $j\neq i$. Thus letting $C_\ell$ act by multiplication of $\zeta_p$ on $V_i$ and trivially on $V_j$ for $j\neq i$ yields a faithful representation of $\tilde{G}:=G\times C_\ell$ on $V$. Now apply the induction hypothesis to $\tilde{G}$. \qed
\end{proof}

\subsection{Minimal dimension of faithful representations}
We define the \emph{representation dimension} of $G$ over $k$ as follows:
\begin{defn}
$\rdim_k G:=\min\{\dim V \mid V \text{ faithful representation of } G \text{ over } k \}.$
\end{defn}
This new numerical invariant gives an upper bound for $\edim_k G$. In certain cases the two invariants of $G$ coincide, e.g. for $p$-groups when $k$ contains a primitive $p$-th root of unity, see \cite[Theorem 4.1]{KM}. \par
\begin{defn} Let $A$ be an abelian subgroup of $G$ and $\chi \in A^\ast:=\Hom(A,\bar{k}^\ast)$.
$$\rep^{(\chi)}(G):=\{V \text{ irreducible representation of } G \mid (V\otimes \bar{k})|_A \supseteq \bar{k}_{\chi} \},$$ where $\bar{k}_{\chi}$ is the one-dimensional representation of $A$ over $\bar{k}$ on which $A$ acts via $\chi$.
\end{defn}
To every group $G$ and field $k$ we associate the following function:
$$f_{G,k} \colon A^\ast \to \NNN_0, \quad \chi \mapsto \min \{\dim V \mid V \in \rep^{(\chi)}(G)\},$$ where $A=\socAb(G)$. \par
From Lemma \ref{le:socrep} we get the following
\begin{cor}
If the socle $C=\soc G$ of $G$ is abelian and $\Chr k \nmid |C|$, then 
$$\rdim_k G = \min \left\{\sum_{i=1}^r f_{G,k}(\chi_i)\right\}$$
taken over all $r \in \NNN$ and all systems of generators $(\chi_1,\dotsc,\chi_r)$ of $C^\ast$ viewed as a $\ZZZ G$-module.
\end{cor}
It may happen that every faithful representation of minimal dimension has more decomposition factors than needed in minimum to create a faithful representation. However in the following situation that doesn't occur and we can describe faithful representations of minimal dimensions more precisely. Recall the definition of a minimal basis introduced in \cite{KM}:
\begin{defn}
Let $C$ be a vector space over some field $F$ of dimension $r\in \NNN_0$ and let $f\colon C \to \NNN_0$ be any function. An $F$-basis $(c_1,\dotsc,c_r)$ of $C$ is called \emph{minimal relative to $f$} if 
\begin{equation}
\label{eq:minsys}
f(c_i) = \min\left\{f(c)\mid c \in C \setminus \langle c_1,\dotsc,c_{i-1}\rangle \right\}, 
\end{equation} for $i=1,\dotsc,r$ where for $i=1$ we use the convention that the span of the empty set is the trivial vector space $\{0\}$. \par
\end{defn}

\begin{prop}
\label{pr:minrep}
Let $G$ be a group whose socle $C:=\soc G$ is a central $p$-subgroup for some prime $p$ and assume $\Chr k \neq p$. Let $V$ be any representation of $G$ and let $V_1,V_2,\dotsc,V_r$ be its irreducible composition factors ordered increasingly by dimension. Choose characters $\chi_1,\dotsc,\chi_r \in C^\ast = \Hom(C,\bar{k}^\ast)$ such that $V_i \in \rep^{(\chi_i)}(G)$. Then $V$ is faithful of dimension $\rdim_k G$ if and only if $r=\rank C$ and $(\chi_1,\dotsc,\chi_r)$ forms a minimal basis of $(C^\ast,f_{G,k})$ with $f_{G,k}(\chi_i)=\dim V_i$. The dimension vector $(\dim V_1,\dotsc,\dim V_r)$ is unique amongst faithful representations of dimension $\rdim_k G$.
\end{prop}
\begin{proof}
Since $p \nmid |C|$ we may replace $V$ by its associated graded representation $V_1\oplus\dotsc \oplus V_r$ without changing faithfulness, decomposition factors and dimension. Thus we will assume that $V$ is completely reducible. \par
First assume that $V$ is faithful and $\rdim_k G = \dim V$. Then the characters $\chi_1,\dotsc,\chi_r$ clearly generate $C^\ast$ and in particular $r\geq\rank C$. Let $j\in \{0,\dotsc,r\}$ be maximal such that $(\chi_1,\dotsc,\chi_r)$ is part of a minimal basis of $C^\ast$. We want to show that $j=r$. Assume to the contrary that $j<r$. Hence there exists $\chi \in C^\ast \setminus \langle \chi_1,\dotsc\chi_j\rangle $ and $W \in \rep^{(\chi)}(G)$ such that $\dim W < \dim V_i$ for all $i>j$. By elementary linear algebra there exists $i>j$ such that $\chi_1,\dotsc,\chi_{i-1},\chi,\chi_{i+1},\dotsc,\chi_r$ generate $C^\ast$ as well. Let $V':=V_1\oplus \dotsb V_{i-1}\oplus W \oplus V_{i+1} \oplus \dotsb \oplus V_r$. Then $\dim V' < \dim V$ and $V'$ is faithful, because $V'$ is faithful restricted to $C$ and every normal subgroup of $G$ intersects $C=\soc(G)$. This contradicts to $\dim V = \rdim_k G$.\par
Now assume that $(\chi_1,\dotsc,\chi_r)$ and $(\chi'_1,\dotsc,\chi'_r)$ form two minimal bases of $C^\ast$. We show that $f_{G,k}(\chi_i)=f_{G,k}(\chi'_i)$ for all $i=1\dotsc r$. Let $j\in \{0,\dotsc, r\}$ be the last index where $(f_{G,k}(\chi_1),\dotsc,f_{G,k}(\chi_j))$ and $(f_{G,k}(\chi'_1),\dotsc,f_{G,k}(\chi'_j))$ coincide. Assume $j<r$ and assume $f_{G,k}(\chi'_{j+1})<f_{G,k}(\chi_{j+1})$. Then $\langle \chi_1,\dotsc, \chi_j \rangle \neq \langle \chi'_1,\dotsc, \chi'_j \rangle$. Hence there exists $s \in \{1,\dotsc,j\}$ such that  $\chi'_s \notin \langle \chi_1,\dotsc, \chi_j \rangle$. Then $f_{G,k}(\chi_{j+1})>f_{G,k}(\chi'_{j+1})\geq f_{G,k}(\chi'_s)$, which contradicts to the definition of minimal basis. This implies uniqueness of the dimension vector and the converse to the above implication.
\qed
\end{proof}
\begin{rem}
Under the assumptions of Proposition \ref{pr:minrep} let $(\chi_1,\dotsc,\chi_r)$ be a minimal basis of $C^\ast$ and $1\leq i_1<i_2<\dotsb<i_m<r$ be the positions of jumps in the vector $(f_{G,k}(\chi_1),\dotsc,f_{G,k}(\chi_r))$, i.e. the indices $i$ where $f_{G,k}(\chi_{i})<f_{G,k}(\chi_{i+1})$. The argument in the proof of Proposition \ref{pr:minrep} shows that the subgroups $\langle \chi_1,\dotsc,\chi_{i_j} \rangle$ for $j=1\dotsc m$ do not depend on the choice of the minimal basis $(\chi_1,\dotsc,\chi_r)$. This yields a canonical filtration $C^\ast=A_{m+1}\supsetneq A_m \supsetneq \dotsc \supsetneq A_1 \supsetneq A_0 = \{e\}$ of $C^\ast$ where $\rank A_j = i_j$ for $j=1,\dotsc,m$. It would be interesting to know whether every basis $(\chi_1,\dotsc,\chi_r)$ of $C^\ast$ respecting this grading of $C^\ast$ is a minimal basis, or equivalently if for all $j=0,\dotsc,m$ and $\chi,\chi' \in A_{j+1}\setminus A_j$ the equality $f_{G,k}(\chi) = f_{G,k}(\chi')$ holds.
\end{rem}
\begin{cor}
Let $p$ be a prime and $G_1, \dotsc,G_n$ be groups. Assume that $\Chr k \neq p$ and $\soc G_i$ is a central $p$-subgroup of $G_i$ for $i=1,\dotsc,n$. Then $$\rdim_k \prod_{i=1}^n G_i = \sum_{i=1}^n \rdim_k G_i.$$ 
\end{cor}
The (statement and the) proof is very similar to \cite[Theorem 5.1]{KM}, which becomes a statement about minimal faithful representations of $p$-groups via \cite[Theorem 4.1]{KM}. Since our situation is more general and we do not require $k$ to contain a primitive $p$-th root of unity, we append the proof.
\begin{proof}
Using induction it suffices to show the case $n=2$. Set $G:=G_1\times G_2$. Taking into account the description of minimal faithful representations of Proposition \ref{pr:minrep} it remains to create a minimal basis $(\chi_1,\dotsc,\chi_{r})$ of $(\soc G)^\ast = (\soc G_1)^\ast \oplus (\soc G_2)^\ast$ for $f_{G,k}$ subject to the condition that each $\chi_i$ is contained in one of $(\soc G_i)^\ast$. Here $r=\rank Z(G) = \rank Z(G_1)+\rank Z(G_2)$. Assume that $(\chi_1,\dotsc,\chi_j)$ is part of a minimal basis such that each $\chi_i$ for $i\leq j$ is contained in one of $(\soc G_i)^\ast$. Choose $\chi \in (\soc G)^\ast \setminus \langle\chi_1,\dotsc,\chi_j\rangle$ with $f_{G,k}(\chi)$ minimal. Decompose $\chi$ as $\chi^{(1)} \oplus \chi^{(2)}$ where $\chi^{(i)} \in (\soc G_i)^\ast$ and choose $W \in \rep^{(\chi)}(G)$ of minimal dimension. The definition of $\rep^{(\chi)}(G)$ means that $\bar{k}_\chi \subseteq W \otimes \bar{k}$. Let $\eps_1$ and $\eps_2$ denote the endomorphism of $G$ sending $(g_1,g_2)$ to $(g_1, e)$ and to $(e,g_2)$, respectively. The representation $\rho_W \circ \eps_i$ contains $\bar{k}_{\chi^({i})}$ and has the same dimension as $W$. Now replace $\chi$ by $\chi^{(i)}$ with $i$ such that $\chi^{(i)}$ lies outside the subgroup of $(\soc G)^\ast$ generated by $\chi_1,\dotsc,\chi_j$. This shows the claim. \qed
\end{proof}
\subsection{Central extensions}
In this subsection we consider central extensions, as investigated in section \ref{sec:CentExt}, from the point of representation theory.
\begin{prop}
\label{part:centextB}
Let $G$ be a semi-faithful group and let $H$ be a central subgroup of $G$ with $H\cap [G,G]=\{e\}$. Let $H'$ be a direct factor of $G/[G,G]$ containing the image of $H$ under the embedding $H\hookrightarrow G/[G,G]$ and assume that $k$ contains a primitive root of unity of order $\exp H'$. Then $$\rdim_k G - \rank Z(G,k) \leq \rdim_k G/H - \rank Z(G/H,k).$$ 
Moreover, if $\soc G$ is a central $p$-subgroup the above inequality is an equality.
\end{prop}
Recall that $G$ is semi-faithful (over $k$) if and only if either $\Chr k = 0$ or $\Chr k = p>0$ and $G$ has no nontrivial normal $p$-subgroups. 
We need some auxiliary results:
\begin{lem}
\label{le:1}
In the situation of the proposition $G/H$ is semi-faithful as well. Moreover there exist characters $\chi_1,\dotsc,\chi_r$ of $G$ such that $\bigcap_{i=1}^r \ker \chi_i \cap H = \{e\}$, where $r=\rank H$. In particular $G$ has a faithful completely reducible representation of the form $V=V'\oplus k^r$ where $V'$ is a completely reducible representation of $G$ with kernel $H$ and $G$ acts on $k^r$ via $g(x_1,\dotsc,x_r) = (\chi_1(g)x_1,\dotsc,\chi_r(g)x_r).$
\end{lem}
\begin{lem}
\label{le:3}
In the situation of the proposition, the quotient homomorphism $\pi \colon G \to G/H$ induces isomorphisms $Z(G)/H\simeq Z(G/H)$ and $Z(G,k)/H\simeq Z(G/H,k)$.
\end{lem}
\begin{proof}[Proof of Lemma \ref{le:1}]
We first show that $G/H$ is semi-faithful over $k$. The case that $\Chr k = 0$ is trivial, hence assume that $k$ has prime characteristic $p$. We now make use of the fact that a group is semi-faithful over $k$ if and only if it does not contain any non-trivial normal abelian $p$-subgroups. Assume that $G/H$ has a normal abelian $p$-subgroup $P \neq \{e\}$. Then the inverse image $B'$ of $P$ under the natural projection is abelian again, since $[B',B'] \subseteq [G,G]\cap H = \{e\}$. Its $p$-Sylow subgroup is then a nontrivial normal abelian $p$-subgroup of $G$. This contradicts to the assumption that $G$ is semi-faithful over $k$. \par
Now let $H'$ be a direct factor of the image of $H$ in $G/[G,G]$ with $\zeta_{\exp H'} \in k$ and let $Z$ denote its complement. Since $k$ contains a primitive root of order $\exp H'$ there exist characters $\tilde{\chi}_1,\dotsc,\tilde{\chi}_r$ of $H'$ such that $\bigcap_{i=1}^r \ker \tilde{\chi}_i$ intersects trivially with the image of $H$ in $H'$. Now define $\chi_i$ by $\chi_i(g)= \chi_i(\pi_2 \pi_1(g))$ where $\pi_1 \colon G \to G/[G,G]$ and $\pi_2 \colon G/[G,G]\simeq H'\times Z \to H'$ are the obvious projection homomorphisms. By construction $\bigcap_{i=1}^r \ker \chi_i \cap H = \{e\}$. \qed
\end{proof}
\begin{rem}
Actually one can show that the conditions of Proposition \ref{part:centextB} are equivalent to the existence of a faithful representation of $G$ of the form given in Lemma \ref{le:1}. The most economical choice for $H'$ is the (unique up to isomorphism) maximal subgroup of $G/[G,G]$ subject to the condition $\soc H' = \soc H$, or in other words, such that for every prime $p$ the $p$-Sylow-subgroup of $H'$ contains the $p$-Sylow-subgroup of $H$ and has the same rank. 
\end{rem}
\begin{proof}[Proof of Lemma \ref{le:3}]
Restricting $\pi$ to $Z(G)$ and $Z(G,k)$ we get homomorphism $Z(G)\to Z(G/H)$ and $Z(G,k)\to Z(G/H,k)$. It remains to show that the two maps are surjective. The map $Z(G)\to Z(G/H)$ is easily seen to be surjective, because if some $g \in G$ commutes with any other $g'\in G$ up to elements of $H$, then it is central, because $[G,G]\cap H=\{e\}$. For the second map let $\pi_2\colon G/[G,G] = Z\times H' \to H'$ denote the projection and consider the homomorphism $G \to G/H\times H', g \mapsto (\pi(g),\pi_2(g[G,G]))$, which is injective. If $\pi(g) \in Z(G/H,k)$ then $k$ contains a primitive root of unity of order $\ord(\pi(g))$ as well as a primitive root of unity of order $\exp H'$. Thus $k$ contains a primitive root of unity of order $\ord g$, whence $g \in Z(G,k)$. \qed
\end{proof}
\begin{proof}[Proof of Proposition \ref{part:centextB}]
Using induction on the order of $H$ we may assume that $H$ is cyclic: The case that $|H|=1$ is clear. If $|H|>1$ then $H$ contains some cyclic subgroup $H_0 \subsetneq H$. By Lemma \ref{le:1} $G/H_0$ is semi-faithful. If $H'$ is a direct factor of $G/[G,G]$ containing $H$ then it contains $H_0$ as well. Moreover $H'/H_0$ is a direct factor of $(G/H_0)/[G/H_0,G/H_0]$ and its exponent is no larger than the exponent of $H'$. Induction yields for the subgroups $H_0 \subseteq G$ and $H/{H_0} \subseteq G/{H_0}$:
\begin{align*}
\rdim_k G - \rank Z(G,k) & \leq \rdim_k G/{H_0} - \rank Z(G/{H_0},k) \\
\rdim_k G/{H_0} - \rank Z(G/{H_0},k) & \leq \rdim_k G/H - \rank Z(G/H,k),
\end{align*}
with equality if $\soc G$ (and therewith $\soc (G/H)$) is a central $p$-subgroup. Combining the two lines shows the claim.
\par We assume now that $H$ is cyclic. 
Let $V$ be a faithful representation of $G/H$ with $\dim V = \rdim_k G/H$. By Lemma \ref{le:arep} we may assume that $V$ is completely reducible, $V=\bigoplus_{i=1}^n V_i$ for some $n \in \NNN$ and irreducible representations $V_i$. We must construct a faithful representation of $G$ of dimension $\dim V + \rank Z(G,k) - \rank Z(G/H,k)$. By (the proof of) Lemma \ref{le:1} there exists a faithful representation of $G$ of the form $V\oplus k_{\chi}$ where $\chi$ is a character whose restriction to $H$ is faithful. \par
If $\rank Z(G,k) = \rank Z(G/H,k) + 1$ this does the job. Otherwise $\rank Z(G,k) = \rank Z(G/H,k)$ and we will consider representations $V_{m_1,\dotsc,m_n} := \bigoplus_{i=1}^n V_i \otimes \chi^{m_i}$ for $m_1,\dotsb,m_n \in \ZZZ$. Clearly $V_{m_1,\dotsc,m_n}$ has the right dimension. We will choose $m_1,\dotsc,m_n$ such that $V_{m_1,\dotsc,m_n}$ becomes faithful. In general let $g$ act trivially on $V_{m_1,\dotsc,m_n}$. This implies that for each $i$ the element $g$ acts like $\chi^{-m_i}$ on $V_i$. In particular the image of $g$ in $G/H$ is an element of $Z(G/H,k)$. Since $Z(G/H,k) \simeq Z(G,k)/H$ under the canonical projection this implies that $g \in Z(G,k)$. Hence $V_{m_1,\dotsc,m_n}$ is a faithful representation of $G$ if and only if it is faithful restricted to $Z(G,k)$. \par
The elements of $Z(G,k)$ act through multiplication with characters $\chi_1,\dotsc,\chi_n$ of $Z(G,k)$ on $V_1,\dotsc,V_n$. Let $\hat{\chi}$ denote the restriction of $\chi$ to $Z(G,k)$. Then the elements of $Z(G,k)$ act through the characters $\chi_1 \hat{\chi}^{m_1},\dotsc,\chi_n \hat{\chi}^{m_n}$ on the irreducible components of $V_{m_1,\dotsc,m_n}$. Using (the second part) of the following Lemma \ref{le:eldiv} we find $m_1,\dotsc,m_n$ such that $\chi_1 \hat{\chi}^{m_1},\dotsc,\chi_n \hat{\chi}^{m_n}$ generate the whole group $\Hom(Z(G,k),\Gm)$ of characters. Then $V_{m_1,\dotsc,m_n}$ is faithful restricted to $Z(G,k)$, hence, as previously observed, faithful for $G$. \par

Now assume that $C:=\soc G$ is a central $p$-group. It then consist precisely of the central elements of exponent $p$ of $G$. We want to show 
$\rdim_k G/H \leq \rdim_k G - (\rank Z(G,k) - \rank Z(G/H,k))$. By assumption $k$ contains a primitive root of unity of order $|H|$ and we may assume $H\neq \{e\}$, hence $\zeta_p \in k$. Let $V=\bigoplus_{i=1}^r V_i$ be a faithful representation of $G$ with $\rdim_k G = \dim V$ and each $V_i$ irreducible. There exist characters $\chi_1,\dotsc,\chi_r \in C^\ast:=\Hom(C,k^\ast)$ such that $cv_i = \chi_i(c)v_i$ for $c \in C$ and $v_i \in V_i$. Faithfulness of $V$ is equivalent to the statement that $\chi_1,\dotsc,\chi_r$ generate $C^\ast$. In particular $r=\rank Z(G) = \rank C$, since $V$ is minimal. Now as in the first part of the proof let $\chi \in \Hom(G,k^\ast)$ be a character which is faithful restricted to $H$. By elementary linear algebra there exists $i \in \{1,\dotsc,r\}$ such that $\chi_1,\dotsc,\chi_{i-1},\chi|_H,\chi_{i+1},\dotsc,\chi_r$ is a basis of $C^\ast$. Replacing $V_i$ by $k_{\chi}$ we get a faithful representation of $G$ of minimal dimension which is of the form $V' \oplus k_{\chi}$. Moreover by multiplying the irreducible components of $V'$ with suitable powers of $\chi$ we may assume that $H$ acts trivially on $V'$. Then the representation $V'':=k_{\chi^|H|} \oplus V'$ is a faithful representation of $G/H$. This establishes the inequality $\rdim_k G/H \leq \rdim_k G - (\rank Z(G) - \rank Z(G/H))$ in case that $\rank Z(G)=\rank Z(G/H)$. In the other case $\rank Z(G)=\rank Z(G/H) + 1$. In that case $\soc(G/H) \simeq C/(H\cap C)$, which is faithfully represented on $V'$, turning $V'$ into a faithful representation of $G/H$ of dimension $\rdim_k G -1$. This finishes the proof. \qed
\end{proof}
\begin{lem}
\label{le:eldiv}
\begin{enumerate}
\item Let $A$ be an abelian group generated by $a_1,\dotsc,a_n \in A$.
Then if $\rank A < n$ there exist $e_1,\dotsc,e_n \in \ZZZ$ co-prime such that $\sum_{i=1}^n e_i a_i = 0$. 
\item Let $A$ be an abelian group generated by elements $c_1,\dotsc,c_n,h$. Assume that $\rank A \leq n$. Then there exist $m_1,\dotsc,m_n \in \ZZZ$ such that $A = \langle c_1 + m_1 h,\dotsc,c_n + m_n h \rangle$.
\end{enumerate}
\end{lem}
\begin{proof}
\begin{enumerate}
\item This follows from the elementary divisor theorem applied to the kernel of the map $\ZZZ^n \twoheadrightarrow A$ sending the $i$-th basis vector of $\ZZZ^n$ to $a_i \in A$. 
\item First assume that the order of $h$ is of the form $p^l$ where $p$ is a prime and $l \in \NNN$. Since $\rank A \leq n$ part (A) shows that there exist $e_1,\dotsc,e_n, e_0 \in \ZZZ$ co-prime such that $\sum_{i=1}^n e_i c_i = e_0 h$. Now if $e_0$ is not divisible by $p$ we get that $h \in \langle c_1,\dotsc,c_n\rangle$ and we can set $m_1=\dotsc=m_n=0$. Otherwise there exists $i \in \{1,\dotsc,n\}$ such that $e_i$ is not divisible by $p$. Then choose $m_i$ such that $e_i m_i \equiv 1- e_0 \pmod {p^l}$ and set $m_j=0$ for $j \neq i$. Then $\sum_{j=1}^n e_j (c_j+m_j h) = (e_0+e_i m_i) h = h$, hence $h \in \langle c_1 + m_1 h,\dotsc,c_n+m_n h\rangle$ and it follows that $A=\langle c_1+m_1 h,\dotsc,c_n+m_n h\rangle$. \par
Now if $h$ is arbitrary we decompose it as $h = \sum_{i=1}^s h_i$ where $h_i$ is of order $p^{l_i}$ for some primes $p_1 < \dotsc < p_s$ and $l_1,\dotsc,l_s \in \NNN$ and apply the just proved statement iteratively to $A_j=\langle c_1,\dotsc,c_n,h_1,\dotsc,h_j\rangle$ for $j=1\dotsc s$ with generators taken from the previous step plus $h_j$. This gives elements $m_{i,j}$, $1\leq i \leq n, 1 \leq j \leq s$ with $A_j=\langle c_1 + \sum_{t=1}^j m_{1,t} h_t,\dotsc, c_n + \sum_{t=1}^j m_{n,t} h_t \rangle$. We have $A=A_s$. The Chinese remainder theorem now implies the claim.
\end{enumerate} \qed
\end{proof}

\begin{cor}
\label{cor:prabB}
Let $G$ and $A$ be groups, where $G$ is semi-faithful and $A$ is abelian. Assume that $k$ contains a primitive root of unity of order $\exp A$. Then 
$$\rdim G\times A - \rank Z(G,k)\times A \leq \rdim G  - \rank Z(G,k)$$
with equality if $\soc G$ is a central $p$-subgroup. 
\end{cor}
\begin{proof}
Apply Proposition \ref{part:centextB} to the central subgroup $\{e\} \times A \subseteq G \times A$. \qed
\end{proof}
\section{Relation of covariant and essential dimension}
\label{sec:EsCov}
The following theorem generalizes \cite[Theorem 3.1]{KLS}, which covers the case $k=\CCC$. 
\begin{theorem}
\label{th:coves}
Let $G$ be non-trivial and semi-faithful. Then $\cdim_k G = \edim_k G$ if and only if $Z(G,k)$ is non-trivial. Otherwise $\cdim_k G = \edim_k G + 1$.
\end{theorem}

\begin{rem}
The theorem does not hold if $\Chr k = p$ and $G$ contains a normal $p$-subgroup. Consider for example an elementary abelian $p$-group, which has essential dimension $1$ by \cite[Proposition 5]{Led}, but covariant dimension $2$, as the following argument shows: It is enough to consider the case $G=\ZZZ/p\ZZZ$. Let $V$ denote the $2$-dimensional representation of $G$ where a generator $g \in G$ acts as $g(s,t)=(s,s+t)$. Suppose that there exists a regular faithful covariant $\phi \colon \AAA(V) \to \AAA(V)$ with $X = \overline{\image \phi}$ of dimension $1$. Then any element $g$ induces an automorphism of order $p$ on the normalization of $X$, which is isomorphic to $\AAA^1$. Since in characteristic $p$ no automorphism of $\AAA^1$ of order $p$ has fixed points we get a contradiction.
\end{rem}
The proof of Theorem \ref{th:coves} remains basically the same as in \cite[section 3]{KLS}. We will append it for convenience. 
\begin{proof}[Proof of Theorem \ref{th:coves}]
Let $Z:=Z(G,k)$ and let $V=\bigoplus_{i=1}^n V_i$ be a faithful representation where each $V_i$ is irreducible. The case when $Z$ is trivial follows from Theorem \ref{th:2}, since $M_\phi$ cannot be the zero-matrix for any regular multihomogeneous covariant $\phi \colon \AAA(V) \to \AAA(V)$. Thus assume that $Z$ is non-trivial. Let $\phi \colon \AAA(V) \dashto \AAA(V)$ be a minimal multihomogeneous covariant. \par 
First assume that there exists a row vector $\beta \in \ZZZ^n$ such that all entries of $\alpha:=\beta M_{\phi}$ are strictly positive. We may assume that $\phi$ is of the form $\phi = \frac{\psi}{f}$ where $f \in k[V]^G$ is multihomogeneous and $\psi \colon \AAA(V) \dashto \AAA(V)$ is a (faithful) regular multihomogeneous covariant. Consider $\tilde{\phi} = (f^{\alpha_1} \phi_1,\dotsc,f^{\alpha_n} \phi_n)$. It is of the form $\gamma \phi$ where $\gamma \in \Hom(X(T_V),\MMMM_G(V))$ is of type I relative to $\phi$. Since $\alpha_j > 0$ for all $j$ the covariant $\tilde{\phi}$ is regular. Lemma \ref{le:dimcov} implies $$\cdim_k G \leq \dim \tilde{\phi} \leq \dim \phi = \edim_k G.$$ \par
We reduce to the case above by post-composing with a covariant as in Example \ref{ex:1}.
Let $g \in Z \setminus \{e\}$ and write $M_\phi=(m_{ij})$. Since $V$ is faithful the element $g$ acts non-trivially on some $V_j$. For such $j$ one of the $m_{ij}$'s must be non-zero. Fix $i_0$ and $j_0$ with $m_{i_0j_0} \neq 0$. Then $\phi_{j_0} \neq 0$ and we can find a homogeneous $h \in k[W_{j_0}]^G$ of degree $\deg h > 0$ such that $h \circ \phi_{j_0} \neq 0$. For any $r \in \ZZZ$ consider the covariant
$$\phi' \colon \AAA(V) \dashto \AAA(V), \quad v \mapsto h^r(\phi_{j_0}(v))\phi(v).$$
Since $h \circ \phi_{j_0} \neq 0$ and $\phi$ is faithful, $\phi'$ is faithful, too. Clearly $\dim \phi' \leq \dim \phi = \edim_k G$. Moreover $\phi'$ is multihomogeneous of degree $M_{\phi'}=(m'_{ij})$ where $m'_{ij} = m_{ij} + r \deg h m_{i j_0}$. For suitable $r \in \ZZZ$ this yields a matrix $M_{\phi'}$ where all $m'_{i_0 j}$ for $j=1\dots n$ are strictly positive. Now for $\beta = e_{i_0}$ the entries of $\alpha = \beta M_{\phi}$ are all strictly positive and we are in the case above. \qed
\end{proof}

\section{The central extension theorem}
\label{sec:CentExt}
As announced in the introduction we shall prove a generalization of the central 
extension theorem.
\begin{theorem}
\label{part:centextA}
Let $G$ be a semi-faithful group. Let $H$ be a central subgroup of $G$ with $H\cap [G,G]=\{e\}$. Let $H'$ be a direct factor of $G/[G,G]$ containing the image of $H$ under the embedding $H\hookrightarrow G/[G,G]$ and assume that $k$ contains a primitive root of unity of order $\exp H'$. Then $$\edim_k G - \rank Z(G,k) = \edim_k G/H - \rank Z(G/H,k).$$ 
\end{theorem}
\begin{rem}
Theorem \ref{part:centextA} generalizes the following results about central extensions: \cite[Theorem 5.3]{BR}, \cite[Theorem 4.5]{Ch}, \cite[Corollary 3.7 and Corollary 4.7]{KLS}, as well as \cite[Theorem 7.1 and Corollary 7.2]{BRV2} and \cite[Lemma 11.2]{BRV}. Chang's version generalizes the result of Buhler and Reichstein to fields of arbitrary characteristic. A closer look reveals that it covers precisely the case of Theorem \ref{part:centextA} when $H$ is cyclic of prime order and maximal amongst cyclic subgroups of $Z(G,k)$. The results of \cite{KLS} do not have these additional assumptions, but they only work for groups $G$ with $\rank Z(G)\leq 2$ and are formulated for the field of complex numbers. Brosnan, Reichstein and Vistoli's Lemma 11.2 from \cite{BRV} gives the inequality $\edim_k G \geq \edim_k G/H$. 
Theorem 7.1 from \cite{BRV2} for fields with $\Chr k \nmid |G|$ extends \cite[Theorem 4.5]{Ch} in the sense that it does not assume any more that $H$ has prime order, but still it makes the assumption that $H$ is maximal amongst central cyclic subgroups of $G$. Corollary 7.2 from \cite{BRV2} is restricted to $p$-groups and it assumes that $H$ is a direct factor of $Z(G)$.
\par
If $G$ is a $p$-group then Theorem \ref{part:centextA} can be deduced from the theorem of Karpenko and Merkurjev about the essential dimension of $p$-groups and Proposition \ref{part:centextB}. 
\end{rem}
\begin{proof}[Proof of Theorem \ref{part:centextA}]
As in the proof of Proposition \ref{part:centextB} we may assume that $H$ is cyclic and there is a faithful representation of $G$ of the form $V \oplus k_{\chi}$ where $\chi$ is faithful on $H$ and $V=\bigoplus_{i=1}^n V_i$ is a completely reducible representation with kernel $H$. We prove the two inequalities of the equation $\edim_k G - \edim_k G/H = \rank Z(G,k) - \rank Z(G/H,k)$ separately:
\par
''$\leq$``: Let $\phi \colon \AAA(V) \dashto \AAA(V)$ be a minimal faithful multihomogeneous covariant of $G/H$. Define a faithful covariant of $G$ via
$$\Phi \colon \AAA(V \oplus k_{\chi}) \dashto \AAA(V \oplus k_{\chi}), \quad (v,t) \mapsto (\phi(v),t).$$
Clearly $\Phi$ is multihomogeneous again of rank $\rank M_{\Phi} = \rank M_{\phi} + 1 = \rank Z(G/H,k) + 1$. Moreover by Theorem \ref{th:2},
$$\edim_k G \leq \dim \Phi - (\rank M_{\Phi} - \rank Z(G,k)) = \edim_k G/H - \rank Z(G/H,k) + \rank Z(G,k).$$
\par
''$\geq$``: Let $\phi \colon \AAA(V \oplus k_{\chi}) \dashto \AAA(V \oplus k_{\chi})$ be a minimal faithful multihomogeneous covariant of $G$. Let $m:=|H|$ and consider the $G$-equivariant regular map $\pi \colon \AAA(V \oplus k_{\chi}) \to \AAA(V\oplus k_{\chi^m})$ defined by sending $(v,t) \mapsto (v,t^m)$. It is a quotient of $\AAA(V\oplus k_{\chi})$ by the action of $H$. The composition $\phi':=\pi \circ \phi \colon \AAA(V \oplus k_{\chi}) \dashto \AAA(V\oplus k_{\chi^m})$ is $H$-invariant, hence we get a commutative diagram: 
$$\xymatrix{ \AAA(V\oplus k_{\chi}) \ar@{-->}[rd]^{\phi'} \ar[d]_{\pi} \ar@{-->}[r]^{\phi} & \AAA(V\oplus k_{\chi}) \ar[d]^{\pi}\\
\AAA(V\oplus k_{\chi^m}) \ar@{-->}[r]^{\bar{\phi}} & \AAA(V\oplus k_{\chi^m})
}
$$
where $\bar{\phi}\colon \AAA(V\oplus k_{\chi^m}) \dashto \AAA(V\oplus k_{\chi^m})$ is a faithful $G/H$-covariant. Since $\pi$ is finite the rational maps $\phi,\phi'$ and $\bar{\phi}$ all have the same dimension $\edim_k G$. Moreover $\phi'$ and $\bar{\phi}$ are multihomogeneous as well. The degree matrix $M_{\phi'}$ is obtained from $M_{\phi}$ by multiplying its last column by $m$ and from $M_{\bar{\phi}}$ by multiplying its last row by $m$. Hence $\rank M_{\phi}=\rank M_{\phi'} = \rank M_{\bar{\phi}}$. Application of Theorem \ref{th:2} yields:
$\edim_k G/H - \rank Z(G/H,k) \leq \dim \bar{\phi} - \rank M_{\bar{\phi}} = \edim_k G - \rank Z(G,k)$.This finishes the proof. \qed
\end{proof}
\begin{cor}
\label{cor:prabA}
Let $G$ and $A$ be groups, where $G$ is semi-faithful and $A$ is abelian. Assume that $k$ contains a primitive root of unity of order $\exp A$. Then 
$$\edim_k G \times A - \rank Z(G,k) \times A = \edim_k G - \rank Z(G,k).$$ 
\end{cor}
\begin{proof}
Apply Theorem \ref{part:centextA} to the central subgroup $\{e\}\times A \subseteq G \times A$. \qed
\end{proof}
\begin{exa}
Consider a group $G_0$ which is generated by a normal subgroup $H$ and an element $g \in G_0 \setminus H$. Let $m:=\ord(g)$ and $n:=\ord(gH)$ be the orders of $g$ in $G$ and in the quotient $G/H$. We form the semi-direct product $G:=C_m \ltimes H$ by letting a generator $c$ of $C_m$ act on $H$ via conjugation by $g$. Consider the surjective homomorphism $$\alpha \colon G=C_m \ltimes H \rightarrow G_0 \text{ given by } \alpha(c) = g \text{ and }\alpha(h) = h \text{ for } h \in H.$$ 
Its kernel is generated by $x:=c^n g^{-n}$, hence cyclic of order $r:=m/n$. The elements $c$ and $g^n$ commute in $G$ and $x$ lies in the center of $G$, since $[x,c] = e$ and $[x,h]=(c^n(g^{-n}hg^n)c^{-n}) h^{-1} = (g^n(g^{-n}hg^n)g^{-n})h^{-1} = e$ for $h \in H$. We obtain a central extension
$$1 \to C_r \to G \to G_0 \to 1.$$
The intersection $[G,G] \cap \langle x \rangle$ is trivial, since $[G,G]\subseteq H$. 
Now let $\pi$ be the set of prime divisors of the order of the abelian socle of $G=C_m \ltimes H$ and assume $\Chr k \notin \pi$. Then $$\edim_k C_m \ltimes H = \edim_k G_0 + \rank Z(C_m \ltimes H,k)- \rank Z(G_0,k).$$
\end{exa}
Another application of the central extension theorem is the following:
\begin{cor}
Let $G$ be a semi-faithful group with faithful completely reducible representation $V$. Let $\phi \colon \AAA(V) \dashto \AAA(V)$ be a minimal faithful multihomogeneous covariant. Assume that $k$ contains a primitive root of unity of order $p$ for some prime $p$. Then the rational map $\pi_V \circ \phi \colon \AAA(V) \dashto \PPP(V)$ has exactly dimension $\dim \phi - \rank Z(G,k)$. 
\end{cor}
\begin{proof}
The inequality $\dim \pi_V \circ \phi \leq \dim \phi - \rank Z(G,k)$ was already shown previously. We use saturation to prove the reversed inequality. We may assume that the rank of $Z(G,k)$ equals the rank of its $p$-Sylow subgroup. By Proposition \ref{pr:satrep} $V$ admits a faithful representation of $\tilde{G}:=G\times C_p^{n-r}$ where $n=\dim T_V$ and $r=\rank Z(G,k)=\rank M_{\phi}$. \par
Corollary \ref{cor:trans} implies the existence of $\gamma \in \Hom(X(T_V),\MMMM_G(V))$ such that $\gamma \phi$ is $D$-equivariant for $D=\identity_{T_V}$. This turns $\tilde{\phi}:=\gamma \phi$ into a faithful (multihomogeneous) covariant for $\tilde{G}$. Corollary \ref{cor:prabA} shows that $\dim \tilde{\phi} \geq \edim_k \tilde{G} = \edim_k G + (n-r)$. Since $\pi_V \circ \tilde{\phi} = \pi_V \circ \phi$ we get $\dim \pi_V \circ \phi = \dim \pi_V \circ \tilde{\phi} \geq \dim \tilde{\phi} - n \geq \dim \phi - r$, showing the claim. \qed
\end{proof}

\section{Subgroups and direct products}
\label{sec:subpr}
\begin{prop}
\label{pr:sub}
Let $H \subseteq G$ be a subgroup. Assume that $G$ has a completely reducible faithful representation which remains completely reducible when restricted to $H$. Then
$$\edim_k G - \rank Z(G,k) \geq \edim_k H - \rank Z(H,k).$$
\end{prop}
\begin{proof}
Let $V=\bigoplus_{i=1}^m V_i$ be a faithful representation of $G$ with each $V_i$ irreducible and completely reducible as a representation of $H$ and let $\phi \colon \AAA(V) \dashto \AAA(V)$ be a minimal faithful covariant which is multihomogeneous. By Theorem \ref{th:2} $\rank M_\phi = \rank Z(G,k)$. Now consider $\phi$ as covariant for $H$. By Proposition \ref{pr:refin} the rank doesn't go down replacing $\phi$ by a multihomogenization $H_{\lambda}(\phi)$ with respect to a refinement into irreducible representations for $H$. Hence again by Theorem \ref{th:2} $\edim_k H - \rank Z(H,k) \leq \dim H_{\lambda}(\phi) - \rank M_{H_{\lambda}(\phi)} \leq \dim \phi - \rank M_\phi = \edim_k G - \rank Z(G,k)$. \qed
\end{proof}
\begin{rem}
There exist pairs $(H,G)$ of a group $G$ with subgroup $H$ such that both $H$ and $G$ are semi-faithful over $k$, but none of the completely reducible faithful representations of $G$ restricts to a completely reducible representation of $H$. We found some examples using the computer algebra system \cite{Magma}, the smallest (in terms of the order of $G$) is a pair of the form $H=S_3$, $G=C_2 \ltimes (C_3\ltimes (C_3\times C_3)) $ in characteristic $2$. Also there are examples in order $72$ with $G=Q_8 \ltimes (C_3\times C_3)$ or $G=C_8 \ltimes (C_3\times C_3)$.
\end{rem}
\begin{prop}
\label{pr:dipr}
Let $G_1$ and $G_2$ be semi-faithful groups. Then $$\edim_k G_1 \times G_2 - \rank Z(G_1\times G_2,k) \leq \edim_k G_1 - \rank Z(G_1,k) + \edim_k G_2 - \rank Z(G_2,k).$$
\end{prop}
\begin{proof}
Let $V=\bigoplus_{i=1}^m V_i$ and $W=\bigoplus_{j=1}^n W_j$ be faithful representations of $G_1$ and $G_2$, respectively, where each $V_i$ and $W_j$ is irreducible. Let $\phi_1 \colon \AAA(V) \dashto \AAA(V)$ and $\phi_2 \colon \AAA(W) \dashto \AAA(W)$ be minimal faithful multihomogeneous covariants for $G_1$ and $G_2$. Then $\rank M_{\phi_1} = \rank Z(G_1,k)$ and $\rank M_{\phi_2}= \rank Z(G_2,k)$ by Theorem \ref{th:2}. The covariant $\phi_1 \times \phi_2 \colon  \AAA(V \oplus W) \dashto \AAA(V \oplus W)$ for $G_1 \times G_2$ is again faithful and multihomogeneous with $\rank M_{\phi} = \rank M_{\phi_1} + \rank M_{\phi_2} = \rank Z(G_1,k) + \rank Z(G_2,k)$. Thus, by Theorem \ref{th:2},
\begin{align*} \edim_k G_1 \times G_2 - \rank Z(G_1 \times G_2,k) & \leq \dim \phi - \rank M_\phi \\& = \dim \phi_1 + \dim \phi_2 - \rank Z(G_1,k) - \rank Z(G_2,k).
\end{align*}
Since $\dim \phi_1 = \edim_k G_1$ and $\dim \phi_2 = \edim_k G_2$ this implies the claim. \qed
\end{proof}
\begin{rem}
We do not know of an example where the inequality in Proposition \ref{pr:dipr} is strict.
\end{rem}

\section{Twisting by torsors}
\label{sec:twist}
Let $V=\bigoplus_{i=1}^m V_i$ be a faithful representation of $G$ where each $V_i$ is irreducible and let $\phi \colon \AAA(V) \dashto \AAA(V)$ be a multihomogeneous covariant of $G$ with $\phi_j \neq 0$ for all $j$. We denote by $\PPP(V):= \PPP(V_1)\times \dotsc \times \PPP(V_m)$ the product of the projective spaces. It is the quotient of a dense open subset of $\AAA(V)$ by the action of $T_V$. We write $\pi_V \colon \AAA(V) \dashto \PPP(V)$ for the corresponding rational quotient map. Since $\phi$ is multihomogeneous there exists a unique rational map $\psi\colon \PPP(V) \dashto \PPP(V)$ making the diagram
$$\xymatrix{\AAA(V) \ar@{-->}[r]^{\phi} \ar@{-->}[d]^{\pi_V}& \AAA(V)\ar@{-->}[d]^{\pi_V} \\
\PPP(V) \ar@{-->}[r]^{\psi} & \PPP(V)
}$$ commute. Let $Z:=Z(G,k)$, which acts trivially on $\PPP(V)$ and let $C\subseteq Z$ be any subgroup. We view $\psi$ as an $H:=G/C$-equivariant rational map. We will twist the map $\psi$ (after scalar extension) by some $H$-torsor to get a rational map between products of Severy-Brauer varieties. We summarize the construction and basic properties of the twist construction, cf. \cite[section 2]{Fl}:
\par
Let $K$ be a field and $H$ be a finite group. A \emph{(right-) $H$-torsor} (over $K$) is a non-empty not necessarily irreducible $K$-variety $E$ equipped with a right action of $H$ such that $H$ acts freely and transitively on $E(K_{\sep})$. The isomorphism classes of $H$-torsors (over $K$) correspond bijectively to the elements of the Galois cohomology set $H^1(K,H)$, where an isomorphism class of an $H$-torsor $E$ corresponds to the class of the cocycle $\alpha = (\alpha_\gamma)_{\gamma \in \Gamma_K}$ defined by $\gamma x = x \alpha_\gamma$, where $x$ is any fixed element of $E$ and $\Gamma_K = \End_K(K_{\sep})$ is the absolute Galois-group of $K$. Every $H$-torsor is of the form $\spec L$ where $L/K$ is a Galois $H$-algebra.
\par
Let $X$ be a quasi-projective $H$-variety over $K$. Let $H$ act on the product $E \times X$ by $h(e,x)=(eh^{-1},hx)$. Then the quotient $(E\times X)/H$ exists in the category of $K$-varieties and will be denoted by ${}^E X$. It is called the \emph{twist of the $H$-variety $X$ by the torsor $E$}.
\par
If $X$ and $Y$ are quasi-projective $H$-varieties and $\psi \colon X \dashto Y$ is a rational map, there exists a canonical rational map ${}^E\psi\colon {}^E X\dashto {}^E Y$. Moreover if $Z$ is another quasi-projective variety and $\psi_1 \colon X \dashto Y$ and $\psi_2 \colon Y \dashto Z$ are composable, then ${}^E \psi_1 \colon {}^E X \dashto {}^E Y$ and ${}^E \psi_2 \colon {}^E Y \dashto {}^E Z$ are composable as well with composition ${}^E (\psi_2 \circ \psi_1)$.
\par
Let $A$ be a central simple $K$-algebra on which $H$ acts on the left by algebra-homomorphisms. Let $E$ be a $H$-torsor corresponding to a Galois $H$-algebra $L/K$. The \emph{twist of $A$ by the torsor $E$}, denoted by ${}^E A$ is defined to be the subalgebra of $H$-invariants of $A\otimes_K L$ where $H$ acts via $h(a \otimes l) = ha \otimes hl$. \par
If $E\simeq H$ is the trivial $H$-torsor then the twist ${}^E X$ (resp. ${}^E A$) is isomorphic to $X$ (resp. $A$). The varieties $X$ and ${}^E X$ (resp. the algebras $A$ and ${}^E A$) become isomorphic over a splitting field $K'/K$ of $E$ (i.e. over a field where $E$ has a $K'$-rational point). 
Let $U$ be a $K$-vector space of dimension $n$. The algebra $\End_K(U)$ carries an action from $\PGL(U)$ via conjugation. Isomorphism classes of central simple $K$-algebras of degree $n$ correspond bijectively to the elements of $H^1(K,\PGL(U))$, via the following assignment: For $T \in H^1(K,\PGL(U))$, represented by a cocycle $\alpha=(\alpha_\gamma)_{\gamma \in \Gamma_K}$, the corresponding central simple algebra is defined to be the sub-algebra of invariants of $\End(U)\otimes_K K_{\sep}$ under the action of $\Gamma_K$ twisted through $\alpha$, defined by $\gamma \cdot_\alpha (\phi\otimes \lambda) = (\alpha_\gamma \phi) \otimes (\gamma\lambda)$ for $\phi \in \End(U)$ and $\lambda \in K_{\sep}$. 
\par
The three different notions of twisting are related as follows:
\begin{lem}
\label{le:comparision}
Let $U$ be a $K$-vector space of dimension $n$. The group $\PGL(U)$ acts on $\PPP(U)$ from the right in the obvious way and on $\End(U)$ via conjugation from the left. Let $\beta \colon H \to \PGL(U)$ be a homomorphism and let $E$ be a $H$-torsor over $K$. Let $H$ act on $\PPP(U)$ and on $\End(U)$ via the homomorphism $\beta$. Then ${}^E\PPP(U) \simeq \SB(A)$ where $A:={}^E \End(U)$. Moreover $A$ is isomorphic to the central simple algebra corresponding to the image of $E$ under the map $H^1(K,H) \stackrel{\beta_{\ast}}\to H^1(K,\PGL(U))$.
\end{lem}
\begin{proof}
The first part is \cite[Lemma 3.1]{Fl}. For the second part, let $E=\spec(L)$ for some Galois $H$-algebra $L$ and fix $\iota \in \Hom(L,K_{\sep}) = E(K_{\sep})$. Then the image of $E$ in $H^1(K,\PGL(U))$ is represented by the cocycle $\alpha = (\beta (h_\gamma))_{\gamma \in \Gamma_K}$ where $h_{\gamma} \in H$ is such that ${}^\gamma \iota = \iota h_\gamma$. In other words $\gamma(\iota(\ell)) = (\iota h_\gamma)(\ell) = \iota(h_\gamma \ell)$ for all $\ell \in L$. Recall that $A={}^E \End(U)$ is the sub-algebra of $H$-invariants of $\End(U)\otimes L$ and the twist $B$ of $\End(U)$ by the cocycle $\alpha$ is the sub-algebra of $\Gamma_K$ invariants of $\End(U)\otimes K_{\sep}$ under the action twisted by the cocycle $\alpha$. Consider the homomorphism of $K$-algebras $\eps:= \identity \otimes \iota \colon \End(U)\otimes L \to \End(U)\otimes K_{\sep}$. It is equivariant in the sense that  $\eps(h_\gamma x) = \gamma \cdot_{\alpha} \eps(x)$ for $x \in \End(U)\otimes L$ and $\gamma \in \Gamma_K$. To see this, we may check it for $x=\phi\otimes \ell$ where $\phi \in \End(U)$ and $\ell \in L$. Then $\eps(h_\gamma x)$ = $h_\gamma\phi \otimes \iota(h_\gamma \ell)$ and $\gamma \cdot_{\alpha} \eps(x) = \gamma \cdot_{\alpha} (\phi\otimes \iota(\ell)) = \beta(h_\gamma) \phi \otimes \gamma (\iota(\ell)) = h_\gamma \phi \otimes \iota(h_\gamma \ell)$. This shows that $\eps(A) \subseteq B$. Since $A$ is simple, the homomorphism $\eps$ maps $A$ injectively into $B$. Counting dimensions yields $\eps(A)=B$. Hence $\eps$ establishes an isomorphism of $K$-algebras between $A$ and $B$, showing the claim. \qed
\end{proof}
We will now apply the twist construction to our particular situation. Let $K/k$ be a field extension and $E$ be an $H$-torsor over $K$. Extending scalars to $K$ we may twist the map $\psi_K$ with $E$ and get a rational map ${}^E\psi_K \colon {}^E\PPP(V_K) \dashto {}^E\PPP(V_K)$.
\begin{lem}
\label{le:twist}
${}^E\PPP(V_K) \simeq \prod_{i=1}^m \SB(A_i)$, where $A_i$ is the twist of $\End_K(V_i \otimes K)$ by the $H$-torsor $E$ and $\End_K(V_i\otimes K)$ carries the conjugation action induced from $G$. Moreover the class of $A_i$ in $\Brauer(K)$ coincides with the image $\beta^{E}(\chi)$ of $E$ under the map $$H^1(K,H)\to H^2(K,C)\stackrel{\chi_{\ast}}\longrightarrow H^2(K,\Gm) = \Brauer(K)$$ where $\chi \in C^\ast$ is the character defined by $g v = \chi(g)v$ for $g \in C$ and $v \in V_i$.
\end{lem}
\begin{proof}
The first claim follows from Lemma \ref{le:comparision}. For the second claim (cf.~\cite[Lemma 4.3]{KM}) consider the commutative diagram 
$$\xymatrix{ H^1(K,H) \ar[r] \ar[d] & H^2(K,C) \ar[d]^{(\chi_i)_{\ast}} \\
H^1(K,\PGL(V_i \otimes K)) \ar[r] & H^2(K,\Gm)
}$$ arising from the following commutative diagram with exact rows:

$$\xymatrix{ 1 \ar[r] & C \ar[r] \ar[d]^{\chi_i}& G \ar[r] \ar[d]^{\rho_{V_i\otimes K}} & H \ar[r]  \ar[d]& 1\\
1 \ar[r] & \Gm \ar[r] & \GL(V_i\otimes K) \ar[r] & \PGL(V_i\otimes K) \ar[r] & 1. \\
}$$ This shows that the image $\beta^{E}(\chi_i)$ of a torsor $E$ over $K$ in $H^2(K,\Gm)$ coincides with the Brauer-class of the central simple algebra corresponding to the image of $E$ in $H^1(K,\PGL(V_i\otimes K))$. By Lemma \ref{le:comparision} this is precisely the twist of $\End(V_K)$ by the $H$-torsor $E$. \qed
\end{proof}

\begin{defn}
Let $X$ and $Y$ be smooth projective varieties. The number $e(X)$ is defined as the least dimension of the closure of the image of a rational map $X\dashto X$. \par
Let $\CCCC$ be a class of field extensions of some field $K$. A generic field of $\CCCC$ is a field $E \in \CCCC$ such that for every $L \in \CCCC$ there exists a $k$-place $E \leadsto L$. The canonical dimension of $\CCCC$ is the least transcendence degree over $K$ of a generic field of $\CCCC$, denoted by $\candim(\CCCC)$ (possibly infinite). \par
If $X$ is a $K$-variety or if $D\subseteq \Brauer(K)$ is a subgroup, the canonical dimension of $X$ (resp. $D$) is defined as the canonical dimension of the class of splitting fields of $X$ (resp. $D$), i.e. the class of field extensions $L/K$, for which $X(L)\neq \emptyset$ (resp. for which $D$ lies in the kernel of the homomorphism $\Brauer(K) \to \Brauer(L)$). It is denoted by $\candim(X)$ (resp. $\candim(D)$).
\end{defn}
\begin{lem}[{\cite[Corollary 4.6]{KM06}}]
\label{le:candim}
Let $X$ be a smooth projective $K$-variety. Then $e(X) = \candim(X)$.
\end{lem}
We only need the inequality $e(X)\geq \candim(X)$ which is established as follows: Let $\psi \colon X \dashto X$ be a rational map with $\dim \psi = e(X)$ and let $Y$ be the closure of the image of $\psi$. One can show that $K(Y)$ is a generic splitting field for $X$. Hence $\candim(X) \leq \tdeg_K K(Y) = \dim \psi = e(X)$.
\begin{lem}
\label{le:edcd}
$$\edim_k G - \rank Z(G,k) \geq e\left({}^E \PPP(V_K) \right) = \candim \left({}^E \PPP(V_K) \right) = \candim (\image \beta^E)$$
\end{lem}

\begin{proof}
Let $\phi \colon \AAA(V)\dashto \AAA(V)$ and $\psi \colon \PPP(V) \dashto \PPP(V)$ be as in the beginning of this section and assume that $\phi$ is minimal, i.e. $\dim \phi = \edim_k G$. By functoriality we have $\dim {}^E \psi_K \leq \dim \psi_K$. Hence $$e\left({}^E \PPP(V_K)\right) \leq \dim {}^E\psi_K \leq \dim \psi_K = \dim \psi.$$
We now show that $\dim \psi \leq \dim \phi - \rank Z(G,k)$. Let $X:=\overline{\image \phi}\subseteq \AAA(V)$. The fibers of $\pi_V|_X \colon X \to \PPP(V)$ are stable under the torus $D_\phi(T_V)\subseteq T_V$. The dimension of $D_\phi(T_V)$ is greater or equal to $\rank Z(G,k)$, since it contains the image of $Z(G,k)$ under $G\hookrightarrow \GL(V)$. Moreover $D_\phi(T_V)$ acts generically freely on $X$. Hence the claim follows by the fiber dimension theorem.
Lemma \ref{le:candim} implies $e\left({}^E \PPP(V_K)\right) = \candim\left({}^E\PPP(V_K)\right)$. The equality $\candim\left({}^E \PPP(V_K)\right) = \candim \image \beta^{E}$ follows easily by Lemma \ref{le:twist}, since it shows that the class of splitting fields of the variety ${}^E \PPP(V_K)$ is identical to the class of common splitting fields of $\beta^E(\chi_1),\dotsc,\beta^E(\chi_m)$. Since $V$ is faithful restricted to $C$ the characters $\chi_1,\dotsc,\chi_m$ generate $C^\ast$. Hence the splitting fields of ${}^E \PPP(V_K)$ are precisely the splitting fields of the image of $\beta^E$ in $\Brauer(K)$. \qed
\end{proof}
\begin{rem}
Lemma \ref{le:edcd} substitutes one part of the proof of the Theorem of Karpenko and Merkurjev about the essential dimension of a $p$-group $G$ when $k$ contains a primitive $p$-th root of unity, saying that $\edim_k G=\rdim_k G$. They show in that case that $\edim_k G \geq \edim [E/G] = \candim (\image \beta^E) + \rank Z(G)$ where $E$ is a generic $G/C$-torsor, $C:=\soc(Z(G))$ and $[E/G]$ is the corresponding quotient stack, see \cite[Theorem 4.2 and Theorem 3.1]{KM}. Our Lemma is more general because $C$ does not need to be a $p$-group. Probably one could also use the stack theoretic approach to show the result of Lemma \ref{le:edcd}, but using multihomogeneous covariants seems more elementary.
\end{rem}
\begin{rem}[The choice of the subgroups $C\subseteq Z(G,k)$]
Karpenko and Merkurjev work with the subgroup of elements of exponent $p$ in $Z(G,k)$. In their setting $G$ is a $p$-group and $\zeta_p \in k$, so $C$ is the smallest subgroup of $Z(G)$ with the same rank as $Z(G)$. In general the best lower bound is obtained with the maximal choice, i.e. with the subgroup $C=Z(G,k)$. This is seen as follows: Set $Z=Z(G,k)$. For a $G/C$-torsor $E'$ over $K$ let $E$ denote its image under $H^1(K,G/C) \to H^1(K,G/Z)$. Then for any $\chi \in Z^\ast$ we have a commutative diagram:
$$\xymatrix{H^1(K,G/C) \ar[r] \ar[d]& H^2(K,C) \ar[r]^{(\chi|_C)_{\ast}} \ar[d] & H^2(K,\Gm) \ar@{=}[d] \ar[r] & \Brauer(K) \ar@{=}[d] \\
H^1(K,G/Z) \ar[r]& H^2(K,Z) \ar[r]^{\chi_\ast} & H^2(K,\Gm) \ar[r] &\Brauer(K)
}$$
Since every element of $C^\ast$ is the restriction of some character $\chi \in Z^\ast$ this shows that $\image (\beta^E) = \image (\beta^{E'})$, hence their canonical dimensions coincide. \par
In general we don't know whether the choice of the subgroup of elements of exponent $p$ in $Z(G,k)$ gives the same lower bound.
\end{rem}

We quote two key results from \cite{KM}:
\begin{theorem}[{\cite[Theorem 2.1 and Remark 2.9]{KM}}]
\label{th:KM1}
Let $p$ be a prime, $K$ be a field and $D\subseteq \Brauer(K)$ be a finite $p$-subgroup of rank $r \in \NNN$. Then $\candim D = \min \left\{\sum_{i=1}^r (\Ind a_i - 1)\right\}$ taken over all generating sets $a_1,\dotsc,a_r$ of $D$. Moreover if $D$ is of exponent $p$ then the minimum is attained for every minimal basis $a_1,\dotsc,a_r$ of $D$ for the function $d \mapsto \Ind d$ on $D$.
\end{theorem}
\begin{theorem}[{\cite[Theorem 4.4 and Remark 4.5]{KM}}]
\label{th:KM2}
Let $1 \to C \to G \to H \to 1$ be an exact sequence of algebraic groups over some field $k$ with $C$ central and diagonalizable. Then there exists a generic $H$-torsor $E$ over some field extension $K/k$ such that for all $\chi \in C^\ast$: $$\Ind \beta^E(\chi) = \gcd \{\dim V \mid V \in \rep^{(\chi)}(G)\}.$$
\end{theorem}
The following corollary works for a slightly larger class of groups than $p$-groups. It becomes \cite[Theorem 4.1]{KM} under the observation that all irreducible representations of $p$-groups have $p$-primary dimension when $\zeta_p \in k$.
\begin{cor}[{cf.~\cite[Theorem 4.1]{KM}}]
\label{cor:edrd}
Let $G$ be an arbitrary group whose socle $C$ is a central $p$-subgroup for some prime $p$ and let $k$ be a field containing a primitive $p$-th root of unity. Assume that for all $\chi \in C^\ast$ the equality $$\gcd \{\dim V \mid V \in \rep^{(\chi)}(G)\} = \min \{\dim V \mid V \in \rep^{(\chi)}(G)\}$$ holds. Then $\edim_k G = \rdim_k G$.
\end{cor}
\begin{proof}
The inequality $\edim_k G \leq \rdim_k G$ is clear. By the assumption on $k$ we have $\rank C = \rank Z(G,k) = \rank Z(G)$. Hence, by Lemma \ref{le:edcd}, it suffices to show $\candim (\image \beta^{E}) = \rdim_k G - \rank C$ for a generic $H:=G/C$-torsor $E$ over a field extension $K$ of $k$. \par
By Theorem \ref{th:KM1} there exists a basis $a_1,\dotsc,a_s$ of $\image \beta^E$ such that $\candim (\image \beta^E) = \sum_{i=1}^s (\Ind a_i - 1)$. Choose a basis $\chi_1,\dotsc,\chi_r$ of $C^\ast$ such that $a_i = \beta^E(\chi_i)$ for $i=1,\dotsc,s$ and $\beta^E(\chi_i)=1$ for $i>s$ and choose $V_i \in \rep^{(\chi_i)}(G)$ of minimal dimension. By assumption $\dim V_i = \gcd \left\{\dim V \mid V \in \rep^{(\chi_i)}(G)\right\}$, which is equal to the index of $\beta^E(\chi_i)$ for the $H$-torsor of Theorem \ref{th:KM2}. \par
Set $V=V_1\oplus \dotsb \oplus V_r$. This is a faithful representation since every normal subgroup of $G$ intersects $C=\soc G$ non-trivially. Then $\candim(\image \beta^E) = \sum_{i=1}^s (\Ind a_i - 1) = \sum_{i=1}^r \Ind \beta^E(\chi_i) - \rank C = \sum_{i=1}^r \dim V_i - \rank C = \dim V - \rank C \geq \rdim_k G - \rank C$.  The claim follows. \qed
\end{proof}
The following was conjectured in case of cyclic subgroups of the Brauer group and proved (over fields of characteristic $0$) for cyclic groups of order $6$ in \cite{CKM}.
\begin{conj}
\label{conj:candim}
Let $D \subseteq \Brauer(K)$ be a finite subgroup. Then $$\candim D = \sum_p \candim D(p),$$ where $D(p)$ denotes the $p$-Sylow subgroup of $D$.
\end{conj}
\begin{rem}
Brosnan, Reichstein and Vistoli asked the following question in \cite[section 7]{BRV}: ``Let X and Y be smooth projective varieties over a field K. Assume that there are no rational functions $X \dashto Y$ or $Y \dashto X$. Then is it true that $e(X \times Y ) = e(X) + e(Y )$?'' 
It remains true in our case that ``a positive answer to this question would imply the conjecture above''.
\end{rem}

\begin{conjcor}
\label{conjcor:1}
Let $G$ be a group whose socle $C:=\soc G$ is central and let $k$ be a field containing a primitive $p$-th root of unity for every prime $p$ dividing $|C|$. Assume that for all $\chi \in C^\ast$ of prime order $\min \dim W = \gcd \dim W$ taken on both sides over all $W \in \rep^{(\chi)}(G)$.
Then $$\edim_k G = \dim V - \sum_p \rank C(p) + \rank C,$$
where $V=\bigoplus V_p$ is a faithful representation of $G$, the direct sum being taken over all primes $p$ dividing $|C|$, and $V_p$ is of minimal dimension amongst representations of $G$ whose restriction to $C(p)$ is faithful. 
\end{conjcor}
\begin{exa}
Using the computer algebra systems \cite{Magma} and \cite{Gap} (and \cite{Sage} to combine the two) we found several examples of non-nilpotent groups for which \cite[Theorem 1.3]{CKM} applies when $k$ is a field containing $\QQQ(\zeta_3)$. These are groups (of order $432$) with $\soc G=Z(G) \simeq C_6$ whose Sylow $2$- and $3$-subgroup have essential dimension $2$ and $3$, respectively. Corollary \ref{conjcor:1} gives for their essential dimension $\edim_k G = (2+3)-2+1=4$. 
\end{exa}
\begin{proof}
``$\leq$'': Consider the multihomogeneous covariant $\identity \colon \AAA(V) \to \AAA(V)$. Theorem \ref{th:2} implies $\edim_k G \leq \dim \identity - (\rank M_{\identity} - \rank Z(G,k)) = \dim V - \sum_p \rank C(p) + \rank C$. \par
``$\geq$'': Choose a generic $G/C$-torsor $E$. Then $\edim_k G \geq \candim (\image \beta^E) + \rank C$, by Lemma \ref{le:edcd}. The $p$-Sylow subgroup of the image of the abelian group $C=\bigoplus_p C(p)$ equals $\beta^E(C(p))$. Conjecture \ref{conj:candim} implies that $\candim \image \beta^E = \sum_p \candim \beta^E(C(p))$, which can be computed with the help of Theorems \ref{th:KM1} and \ref{th:KM2}. Similarly as in the proof of Corollary \ref{cor:edrd} we get the claim, using the replacement of $\gcd$ by $\min$. \qed
\end{proof}
\begin{exa}
Let $G$ be nilpotent, i.e. the direct product of its Sylow subgroups $G(p)$, $p$ prime. Assume that $k$ contains a primitive $p$-th root of unity for every prime $p$ dividing $|G|$. Then Conjecture \ref{conj:candim} and its corollary imply 
$$\edim_k G = \sum_p (\rdim_k G(p) - \rank C(p)) + \rank C.$$
\end{exa}
\section{Normal elementary $p$-subgroups}
\label{sec:FunctorialEd}
Suppose that we are in the case of a non-semi-faithful group $G$. Recall that this happens precisely when $\Chr k = p>0$ and $G$ contains a nontrivial normal $p$-subgroup $A$. Replacing $A$ by the elements of $Z(A)$ of exponent $p$ (which is again normal in $G$) we may assume that $A$ is $p$-elementary. In particular $\edim_k A = 1$ by \cite[Proposition 5]{Led}. We would like to relate $\edim_k G$ and $\edim_k G/A$ and use this iteratively to pass to the semi-faithful case. \par
Merkurjev's description of essential dimension as the essential dimension of the Galois cohomology functor $H^1(\_,G)$ from the category of field extensions of $k$ to the category of sets (see \cite{Favi}) gives the following:
\begin{prop}
\label{pr:estimate}
If $A$ is an elementary $p$-group contained in the center of $G$ and if $\Chr k = p$ then 
\begin{equation*}
\label{eq:edrel}
\edim_k G/A \leq \edim_k G \leq \edim_k G/A+1. \tag{$\ast$}
\end{equation*} 
\end{prop}
\begin{proof}
Since $A$ is central there is the following exact sequence in Galois cohomology:
$$1 \to H^1(\_,A) \to H^1(\_,G) \to H^1(\_,G/A) \to H^2(\_,A)=1.$$
Thus $H^1(\_,G) \to H^1(\_,G/A)$ is a surjection of functors. In particular $\edim_k G/A \leq \edim_k G$ by \cite[Lemma 1.9]{Favi}. \par 
We have an action of $H^1(\_,A)$ on $H^1(\_,G)$ as follows: Let $K/k$ be a field extension and let $[\alpha] \in H^1(K,A)$ and $[\beta] \in H^1(K,G)$ and set $[\alpha]\cdot [\beta] :=[\alpha \beta] \in H^1(K,G)$. Since $A$ is a central $\alpha \beta$ satisfies the cocyle condition and its class in $H^1(K,G)$ does not depend on the choice of $\alpha$ and $\beta$. Moreover it is well known that two elements of $H^1(K,G)$ have the same image in $H^1(K,G/A)$ if and only if one is transformed from the other by an element of $H^1(K,A)$, see \cite{Se}. Thus we have a transitive action on the fibers of $H^1(K,G)\twoheadrightarrow H^1(K,G/A)$, and this action is natural in $K$. That means we have a fibration of functors $$H^1(\_,A) \leadsto H^1(\_,G) \twoheadrightarrow H^1(\_,G/A).$$ Now \cite[Proposition 1.13 ]{Favi} yields $\edim_k G\leq \edim_k G/A+\edim_k A = \edim_k G/A+1$. \qed
\end{proof}
\begin{rem}
If $G$ is a $p$-group and $A$ is a (not necessarily central) elementary abelian $p$-subgroup contained in the Frattini subgroup of $G$ then \cite{Led2} gives the relations \eqref{eq:edrel} as well.
\end{rem}
\begin{exa}
Let $G$ denote the perfect group of order $8!=40320$ which is a central extension of $A_8$ by $C_2$. The socle of this group $\soc G = C_2$ is central. \\
{\bf Claim:} 
$\edim_k G = 8$ if $\Chr k \neq 2$ and $\edim_k G \in \{2,3,4\}$ if $\Chr k = 2$.
\begin{proof}
First consider the case when $\Chr k \neq 2$. There exists a faithful irreducible representation of $G$ of degree $8$ with entries in $\mu_2(k)\simeq C_2$. This implies in particular that $\edim_k G \leq 8$. Moreover one may check using a Computer algebra system like \cite{Magma} or \cite{Gap} that the degree of every faithful irreducible representation of $G$ is a multiple of $8$. The faithful irreducible representations of $G$ are precisely the elements of $\rep^{(\chi)}(G)$ where $\chi$ is the non-trivial character of $\soc G = C_2$. Hence the claim follows with Corollary \ref{cor:edrd}. \par
Now consider the case of $\Chr k = 2$. Proposition \ref{pr:estimate} implies that $\edim_k A_8 \leq \edim_k G \leq \edim_k A_8+1$. The essential dimension of $A_8 \simeq \GL_4(\FFF_2)$ is either $2$ or $3$, see \cite[Lemma 5.5 and Theorem 5.6]{Ch}, and the claim follows. \qed
\end{proof}
\end{exa}

\section*{Acknowledgments}
I am grateful to my PhD adviser Hanspeter Kraft for many fruitful discussions about covariant and essential dimension and the content of this paper.

\end{document}